\documentclass[reqno]{amsart}
\usepackage{amsfonts}
\usepackage{amssymb}
\usepackage{amsthm}
\usepackage{upref}
\usepackage{enumerate}

\usepackage{amscd}

\makeatletter
\def\LaTeX{\leavevmode L\raise.42ex
    \hbox{\kern-.3em\size{\sf@size}{0pt}\selectfont A}\kern-.15em\TeX}
 \makeatother
 
\DeclareMathOperator{\clos}{clos}

\numberwithin{equation}{section}

\newtheorem{lemma}{Lemma}[section]
\newtheorem{theorem}[lemma]{Theorem} 
\newtheorem{corollary}[lemma]{Corollary}
\newtheorem{proposition}[lemma]{Proposition}

\theoremstyle{definition}

\newtheorem{example}[lemma]{Example}

\newtheorem{remark}[lemma]{Remark}

\newcommand{\arccosh}{\operatorname{arccosh}}

  \newcommand{\e}{\eqref}

\newcommand{\q}{\quad}

\newcommand{\ov}{\overline}
\newcommand{\wt}{\widetilde}
\newcommand{\z}{\zeta}

\renewcommand{\d}{\delta}

   \newcommand{\sgn}{\operatorname{sgn}}

\renewcommand\Im{\operatorname{Im}}
\renewcommand\Re{\operatorname{Re}}

\newenvironment{pf}{\begin{proof}}{\end{proof}}

\def\qqq{\mathrel{\subset\mkern-15mu\lower.38ex\hbox{${\scriptscriptstyle\rightarrow}$}}}

\let\cal\mathcal

\let\Bbb\mathbb
 
\begin{document}
\title 
[Asymptotic  behavior of orthogonal polynomials]
{Asymptotic  behavior of orthogonal polynomials \\ without the Carleman condition} 
\author{ D. R. Yafaev  }
\address{   Univ  Rennes, CNRS, IRMAR-UMR 6625, F-35000
    Rennes, France and SPGU, Univ. Nab. 7/9, Saint Petersburg, 199034 Russia}
\email{yafaev@univ-rennes1.fr}
 

   \date{28 November 2019}
\thanks {Supported by  project   Russian Science Foundation   17-11-01126}

\begin{abstract}
Our goal  is to find  an asymptotic behavior as $n\to\infty$ of orthogonal polynomials $P_{n}(z)$  defined by the Jacobi recurrence coefficients $a_{n}, b_{n}$. We suppose that the off-diagonal coefficients  $a_{n}$ grow so rapidly
  that the series $\sum  a_{n}^{-1}$ converges, that is, the Carleman condition is violated. With respect to diagonal coefficients  $b_{n}$ we assume that $-b_{n} (a_{n}a_{n-1})^{-1/2}\to 2\beta_{\infty}$ for some $\beta_{\infty}\neq \pm 1$.
  The asymptotic formulas obtained  for $P_{n}(z)$   are quite different from the case $\sum  a_{n}^{-1}=\infty$ when the Carleman condition is satisfied. In particular,  if $\sum  a_{n}^{-1}<\infty$, then the phase factors in these formulas do not depend on the spectral parameter  $z\in{\Bbb C}$.   The asymptotic formulas obtained in the cases $|\beta_{\infty}|<1 $ and $|\beta_{\infty}|>1 $ are  also qualitatively different from each other.  As an application of these results, we find necessary and sufficient conditions for the
   essential self-adjointness of the corresponding minimal Jacobi operator.
      \end{abstract}

  %These results imply, in  particular, that the corresponding Jacobi operator has deficiency indices $(1,1)$ in the first case. In the second case, we find necessary and sufficient conditions for the
%   essential self-adjointness of the minimal Jacobi operator.

%\thispagestyle{empty}

 \maketitle

%***********************************************************
\section{Introduction}
%*

\subsection{Overview}

 Orthogonal polynomials  $P_{n } (z)$ can be defined by  a recurrence relation
   \begin{equation}
 a_{n-1} P_{n-1} (z) +b_{n} P_{n } (z) + a_{n} P_{n+1} (z)= z P_n (z), \q n\in{\Bbb Z}_{+}, \q z\in{\Bbb C}, 
\label{eq:RR}\end{equation}
with the boundary conditions $P_{-1 } (z) =0$, $P_0 (z) =1$. We always suppose that $a_n >0$, $b_n=\bar{b}_n $. Determining $P_{n } (z)$, $n=1,2, \ldots$, successively from \e{eq:RR}, we see that $P_{n } (z)$ is a polynomial of degree $n$: $P_{n } (z)=\gamma_{n}z^n+\cdots$ where $\gamma_{n}= (a_{0}a_{1}\cdots a_{n-1})^{-1}$. 

We are interested in the asymptotic behavior of the polynomials  $P_{n } (z)$ as $n\to\infty$. This is a classical problem investigated under various assumptions including usually the   condition
    \begin{equation}
\sum_{n=0}^\infty a_{n}^{-1}=\infty
\label{eq:Carl}\end{equation}
introduced by T.~Carleman in his book \cite{Carleman}. We are aware of only one paper \cite{Sw-Tr} where the asymptotics of    $P_{n } (z)$ was studied without assumption \e{eq:Carl}; this paper is discussed at the end of this subsection. On the contrary, under assumption \e{eq:Carl}
there is an enormous literature on this subject.
% reviewed very briefly in Sect.~3.3 (?). 
Here we mention   Nevai's approach (see his book \cite{Nev}) which allowed the authors of \cite{Mate} to treat the case $a_{n} \to a_{\infty}>0$, $b_{n} \to 0$ as $n\to\infty$  in such a way   that
\begin{equation}
\sum_{n=0}^\infty \big(| a_{n+1} - a_{n} |+|b_{n+1} - b_{n}|\big)<\infty .
\label{eq:LR}\end{equation}

The  asymptotics  of the polynomials  $P_{n } (z)$ for real $z$ in case of the coefficients $a_{n} \to\infty$
but satisfying the  Carleman condition \e{eq:Carl} was studied in the papers
\cite{Janas} (see Theorem~3.2) and \cite{Apt} (see Theorem~3). The authors of \cite{Janas} solved equations \e{eq:RR}
for $P_{n} (z)$ successively for $n=1,2,\ldots$ (the transfer matrix method) which yielded a representation for $(P_{n} (z),
P_{n+1} (z))$ as a product of $n$ $2\times 2$-matrices. The conditions on the coefficients were rather restrictive in   \cite{Janas}; in particular, it was assumed there that $b_{n}=0$  for all $n$. Under broader assumptions this problem was considered in  \cite{Apt} where Nevai's method was used. 

Under very general assumptions on the coefficients $a_{n}, b_{n}$ the asymptotics  of the polynomials  $P_{n } (z)$ was studied in \cite{Sw-Tr} where however it was assumed that $z\in {\Bbb R}$ and that the coefficients $b_{n}$ are small compared to $a_{n}$.   The Carleman and non-Carleman cases were treated in \cite{Sw-Tr} at an equal footing so that the difference in the corresponding asymptotic formulas for $P_{n } (z)$ was not quite visible. Among the results of the  present paper,
 Corollary~\ref{GPSc} seems to be the closest to the main result, Theorem~C, of  \cite{Sw-Tr}. 

I thank G.~\'{S}widerski for useful discussions.

%Later,  Nevai's method was used in  the paper \cite{Apt} to find asymptotics  of the polynomials  $P_{n } (z)$ in case of the coefficients $a_{n} \to\infty$but satisfying the  Carleman condition \e{eq:Carl}. 

\subsection{Main results}

Our goal is to  find asymptotic formulas for $P_{n } (z)$ as $n\to\infty$ without the  Carleman condition \e{eq:Carl}, that is, in the case
 \begin{equation}
  \sum_{n =0}^{\infty} a_{n}^{-1} <\infty ;
 \label{eq:nc}\end{equation} 
 then the coefficients $a_{n}\to \infty$ faster than $n$.
 Astonishingly, the asymptotics of the orthogonal polynomials in this a priori highly singular case is particularly simple and general. 
  Let us briefly describe some of our main results omitting minor technical assumptions. In addition to \e{eq:nc}, suppose that there exists a finite limit
 \begin{equation}
   -\frac{b_{n}}{2\sqrt{a_{n-1}a_{n}} }=:\beta_{n}\to  \beta_{\infty}  \q  \mbox{where}\q |\beta_{\infty}|\neq 1
\label{eq:Gr}\end{equation}
   as $n\to\infty$. %Without loss of generality, we suppose that $\beta\geq 0$. 
 We distinguish two cases:  $|\beta_{\infty}|<1$ and $|\beta_{\infty}|>1$.
 
 If $|\beta_{\infty}|<1$, we  set 
 \begin{equation}
 \phi_{n}= \sum_{m=0}^{n-1} \arccos\beta_{m}
\label{eq:Grf}\end{equation}
where the sum is restricted to $m$ such that $|\beta_{m}|\leq 1$. For  example,   $\phi_{n}=\pi n /2$ if $b_{n}=0 $ for all $n\in{\Bbb Z}_{+}$.
We show in Theorem~\ref{LC} that, for an arbitrary $z\in{\Bbb C}$, all solutions $F_{n} (z)$ of the Jacobi equation  
    \begin{equation}
 a_{n-1} F_{n-1}  (z) +b_{n} F_{n} (z) + a_{n} F_{n+1} (z)= z F_{n} (z), \q n\in{\Bbb Z}_{+}, 
\label{eq:Jy}\end{equation}
have asymptotic  behavior  
 \begin{equation}
F_{n}(  z )= a_{n} ^{-1/2}   \Big(k_{+}   e^{-i\phi_{n} }+ k_{-}   e^{ i\phi_{n}}+ o( 1)\Big), \q n\to\infty,
\label{eq:A2P}\end{equation} 
  with some constants $k_{\pm}  \in{\Bbb C}$ (depending of course on the solution $\{F_{n}  (z)\}$). The equalities  
$  k_{+}=k_{-}=0$ are possible for the trvial solution $F_{n} (z)=0$ only. We emphasize that  due to condition  \e{eq:nc}, the inclusion  $\{ F_{n}(  z )\}\in \ell^2 ( {\Bbb Z}_{+})$ holds for all $z\in {\Bbb C}$.
  Conversely, for all $k_{\pm}  \in{\Bbb C}$, there exists a solution $F_{n} (z)$ of equation \e{eq:Jy} with asymptotics \e{eq:A2P}. In particular, formula \e{eq:A2P} with some constants $k_{\pm} (z)$  is true for the orthogonal polynomials $P_{n} (z)$.

%  Actually, the same asymptotics  is true for 
 %The constants $k_{\pm} (z)$ in the asymptotic formulas for solutions of   \e{eq:Jy} are of course different for different $F_{n} (z)$. 
  
  %We emphasize that these asymptotics are the same for real $z$ and for $z$ with $\Im z\neq 0$.
  
In the case $|\beta_{\infty}|>1$ we   set 
 \begin{equation}
\varphi_{n}= \sum_{m=0}^{n-1} \arccosh |\beta_{m}|
\label{eq:Grf1}\end{equation}
where the sum is restricted to $m$ such that $|\beta_{m}| \geq 1$.  
Then for all $z\in {\Bbb C}$, we have (see Theorem~\ref{GE1}) an asymptotic relation 
 \begin{equation}
    P_{n}(  z )= {\sf k}(z) a_{n} ^{-1/2} (\sgn \beta_\infty)^n e^{\varphi_{n}}  (1+ o(1)).
\label{eq:A2P1}\end{equation} 
If the Jacobi operator with the coefficients $a_{n}, b_{n}$ is essentially self-adjoint, then 
the coefficient ${\sf k} (z)\neq 0$ unless $z$ is  an eigenvalue of this operator. 
%Thus, $    P_{n}(  z )\to\infty$ exponentially as $n\to\infty$.

Note that according to \e{eq:A2P} and \e{eq:A2P1} the asymptotic behavior of the polynomials $P_{n} (z)$
is the same for all $z\in {\Bbb C}$, both for real $z$ and for $z$ with $\Im z\neq 0$; only the coefficients $k_{\pm}  (z) $ and ${\sf k}  (z) $
depend on $z$. 
%The condition \e{eq:Gr} is very essential for our construction. This will be briefly discussed in Sect.~3.3. \boldsymbol{\tau}

We emphasize that the asymptotic formulas   for the cases  \e{eq:Carl}  and  \e{eq:nc}  are qualitatively different from each other. This  is discussed   in Sect.~5.2 and 5.3.
% where it is assumed that condition \e{eq:Gr} is satisfied with $ |\beta_{\infty}|<1$. 
 
 \subsection{Scheme of the approach}
 
 We are motivated by an analogy of the difference \e{eq:Jy} and differential 
 \begin{equation}
 - (a(x) f ' (x, z) )'+ b(x) f  (x, z)= z f  (x, z), \q  x>0,
\label{eq:Schr}\end{equation}
equations where  $a(x) > 0$  and $b(x)$ is real. To a large extent, $x$, $a(x)$ and $b(x)$ here play the roles of the parameters $n$, $a_n$ and $b_{n}$ in the Jacobi equation \e{eq:Jy}. In the case $a(x)=1$, $b\in L^1 ({\Bbb R}_{+})$,  
equation \e{eq:Schr} has a solution $f(x,z)$,  known as the Jost solution, behaving like $e^{i\sqrt{z}  x}$, $\Im \sqrt{z}\geq 0$, as $x\to\infty$. For rather general coefficients $a(x)$,  $b(x)$,   equation \e{eq:Schr} has solutions $f(x,z)$ with asymptotics given by the classical Liouville-Green formula (see   Chapter~6 of the book \cite{Olver}).  In the case $a(x)\to a_{\infty}>0$,  $b(x)\to 0$,  the Liouville-Green formula was  simplified in \cite{Y-LR} which yields solutions $ f  (x, z)$ of \e{eq:Schr} with asymptotics
 \begin{equation}
 f  (x, z)\sim    \exp  \Big(-\int_{0} ^x \big(\frac{b(y)-z}{a(y)}\big)^{1/2} dy\Big)=: Q  (x, z),\q
\Re \big(\frac{b(y)-z}{a(y)}\big)^{1/2} \geq 0 ,
\label{eq:Ans}\end{equation}
as $x\to\infty$. Note that the function $ Q  (x, z)$ (the Ansatz for the  Jost solution $f  (x, z)$) satisfies equation \e{eq:Schr} with a sufficiently good accuracy. Formula \e{eq:Ans} was modified in \cite{JLR} for Jacobi equations \e{eq:Jy} where $a_{n}\to a_{\infty}>0$, $ b_{n} \to 0$ as $n\to\infty$ and condition \e{eq:LR} is satisfied. This permitted to find asymptotics of the orthogonal polynomials $P_{n} (z)$ for such coefficients $a_{n} , b_{n}$ in a very natural way.

We are applying the same scheme in the non-Carleman case.   Let us briefly discuss the main steps of our approach. 

A. 
First, we forget about the orthogonal polynomials $P_{n} (z)$
and distinguish solutions (the Jost solutions) $f_{n}  (z)$ of the  difference  equation \e{eq:Jy} 
 by their asymptotics as $n\to\infty$. This requires a construction of an Ansatz $Q_{n} (z)$   for   the Jost solutions.

  %For differential equation \e{eq:Schr} where $a(x)\to a_{0}>0$  and $b(x )\to 0$  as $x\to\infty$, it was reviewed recently in \cite{Y-LR}.
 
  B. 
%The first step is to find of an Ansatz $q_{n} (z)$   for   the Jost solutions of equation \e{eq:Jy}. 
Under assumption \e{eq:nc} this construction (see Sect.~2.3) is very explicit and, in particular, does not depend on $z\in{\Bbb C}$. In  the case  $|\beta_{\infty} |< 1$, we set
  \begin{equation}
Q_{n}= a_{n}^{-1/2}e^{-i\phi_{n}} 
\label{eq:Ans1}\end{equation}
with the phase $\phi_{n}$ defined by formula \e{eq:Grf}. In  the case  $|\beta_{\infty} |> 1$, the Ansatz equals 
 \begin{equation}
Q_{n}= a_{n}^{-1/2} (\sgn\beta_{\infty})^n e^{-\varphi_{n}}
\label{eq:Ans2}\end{equation}
where the phase $\varphi_{n}$ is given by   \e{eq:Grf1}. It is shown in Sect.~2.4 that in both cases
  the relative remainder
\begin{equation}
 r_{n} (z)  : =( \sqrt{ a_{n-1} a_{n}} Q_{n}  )^{-1} \big(a_{n-1} Q_{n-1}   + (b_{n}-z)Q_{n}   + a_{n} Q_{n+1}  \big), \q n\in{\Bbb Z}_{+},
\label{eq:Grr}\end{equation}
%{eq:re2}
belongs to $ \ell^1 ({\Bbb Z}_{+})$.   
  At an intuitive level, the fact that the Ans\"atzen \e{eq:Ans1} and \e{eq:Ans2} do not depend on $z\in{\Bbb C}$ can be explained by the fast growth of the coefficients $a_{n}$ which makes the spectral parameter $z$  negligible; see Remark~\ref{const}.
  % Sect.~3.3 for more details. 
  
  Actually,  the Ans\"atzen we use (especially, the amplitude factor $a_{n}^{-1/2} $)  are only distantly similar to the Liouville-Green Ansatz for the Schr\"odinger equation  \e{eq:Schr}.  

C.
Then we make in Sect.~2.5  a multiplicative change of variables
 \begin{equation}
   f_{n} (z)= Q_{n}  u_{n} (z ) 
      \label{eq:Jost}\end{equation} 
which permits us to reduce   the difference equation \e{eq:Jy} for  $  f_{n} (z)$  to a Volterra ``integral" equation for the sequence $u_{n} (z)$.   
This equation  depends of course on parameters $a_{n}$,$b_{n}$. In particular, it is   somewhat different in the cases $|\beta_{\infty} |< 1$ and  $|\beta_{\infty} |> 1$. However in both cases this Volterra equation is standardly solved    by iterations in Sect.~3.1 and 3.2. This allows us to prove  that $ u_{n} (z)\to 1$ as $n\to\infty$, and  $ u_{n} (z)$  are analytic functions of $z\in{\Bbb C} $. According to   \e{eq:Ans1} or \e{eq:Ans2} and \e{eq:Jost}   this yields (see Sect.~3.3) asymptotics of  the Jost solutions $ f_{n} (z)$.

% All standard statements about spectral properties of the Jacobi operator $J$ and an asymptotic behavior of the polynomials $P_{n} (z)$ are easy consequences of   this analytic result.
   
   D. 
   The sequence 
   \[
   \tilde{f}_{n} (z )=\ov{f_{n}(\bar{z})}
   \]
also    satisfies  the equation \e{eq:Jy}. In the case $|\beta_{\infty} |< 1$, the solutions  $ f_{n} (z )$ and $   \tilde{f}_{n} (z )$ are linearly independent. Therefore it follows from \e{eq:Ans1}   that all solutions  of the Jacobi equation \e{eq:Jy} have  asymptotic behavior \e{eq:A2P}.
   
   In the case $|\beta_{\infty}| > 1$,   a solution $ g_{n} (z ) $ of \e{eq:Jy} linearly independent with $ f_{n} (z ) $ can be constructed by an explicit formula
    \begin{equation} 
g_{n} (z ) = f_{n} (z )\sum_{m=n_{0}}^n (a_{m-1} f_{m-1}(z) f_{m}(z))^{-1},\q n\geq n_{0} ,
\label{eq:GEg}\end{equation}
where $n_{0}=n_{0}(z)$ is a sufficiently large number. This solution grows exponentially as $n\to\infty$,
 \[
   g_{n} (  z )=  \frac{1} {2 \sqrt{\beta_{\infty}^2-1}} a_{n} ^{1/2} (\sgn\beta_{\infty})^{n+1} e^{\varphi_{n}} (1+ o(1)) .
\]
%label{eq:A2P3}\end{equation}
Since $g_{n} (  z )$ is linearly independent with $f_{n} (  z )$, the polynomials $P_{n} (z) $  are linear combinations of $ f_{n} (  z )$ and $g_{n} (  z )$ which leads to the formula \e{eq:A2P1}.  Note that  \e{eq:GEg}
is a discrete analogue of formula (1.36) in Chapter~4 of the book \cite{YA} for the Schr\"odinger equation.

 The scheme described briefly above seems to be quite different from  \cite{Sw-Tr} where the first step was a study of Tur\'an determinants $P_{n}(z)^2 -P_{n-1}(z) P_{n+1}(z)$.
   
%The same approach works also if the Carleman condition \e{eq:Carl} and assumption \e{eq:Gr} are satisfied. In this case, the asymptotic behavior of $P_{n}  (z)$ is quite different from \e{eq:A2P}. It is briefly described in Sect.~5.2, but the proofs will be given elsewhere.
 
 \subsection{Jacobi operators}
 
It is natural (see the book \cite{AKH}) to associate with the coefficients $a_{n} , b_{n} $  
 a  three-diagonal matrix  
  \begin{equation}
{\cal J} = 
\begin{pmatrix}
 b_{0}&a_{0}& 0&0&0&\cdots \\
 a_{0}&b_{1}&a_{1}&0&0&\cdots \\
  0&a_{1}&b_{2}&a_{2}&0&\cdots \\
  0&0&a_{2}&b_{3}&a_{3}&\cdots \\
  \vdots&\vdots&\vdots&\ddots&\ddots&\ddots
\end{pmatrix} 
\label{eq:ZP+}\end{equation}
known as the Jacobi matrix.
Then equation \e{eq:RR} with the boundary condition $P_{-1 } (z) =0$ is
 equivalent to the equation ${\cal J} P( z )= z P(z)$ for the vector $P(z)=\{ P_{n} (z)\}_{n=0}^\infty$. Thus $P(z)$ is the ``eigenvector" of the matrix $\cal J$ corresponding to the ``eigenvalue" $z$. 

Let us now consider Jacobi operators  defined   by matrix \e{eq:ZP+} in the canonical basis $e_{0}, e_{1}, \ldots$ of the space $\ell^2 ({\Bbb Z}_{+})$. 
The minimal Jacobi operator $J_{0}$ is  defined   on a set $\cal D \subset \ell^2 ({\Bbb Z}_{+})$ of vectors with a  finite number of non-zero components.  It is symmetric  in the space $\ell^2 ({\Bbb Z}_{+})$, and its adjoint operator $J^*_{0}$ is given by the same matrix \e{eq:ZP+} on all vectors $f=\{f_{n}\}\in \ell^2 ({\Bbb Z}_{+})$ such that ${\cal J}f \in \ell^2 ({\Bbb Z}_{+})$.
The deficiency indices of the operator $J_{0}$ are either $(0,0)$ (the limit point case) or $(1,1)$ (the limit circle case). If the Carleman condition \e{eq:Carl} holds, then for all $b_{n}$
the operator $J_{0}$ is essentially
 self-adjoint on $\cal D$ so that  $J_{0}$ has the unique self-adjoint extension $J=\clos J_{0}$ (the closure of $J_{0}$).  %We will see that the same is true if
 % assumptions \e{eq:nc} and  \e{eq:Gr}
% with  $|\beta_{\infty}|>1$ are  satisfied.

On the contrary, 
under assumption \e{eq:nc} the deficiency indices of  $J_{0}$ depend on the value of $|\beta_{\infty}|$. If $|\beta_{\infty}|< 1$,  then it follows from \e{eq:A2P} that all solutions of equation \e{eq:Jy} are in $\ell^2 ({\Bbb Z}_{+})$, and 
hence the deficiency indices of the operator $J_{0}$ are $(1,1)$. In this case
  the operator $J_{0}$
   has a one-parameter family of self-adjoint extensions $J\subset J_{0}^*$.  
 Their domains  can be   described explicitly (see Sect.~6.5 of \cite{Schm} or \S 2 of \cite{Simon}) in terms of the orthogonal polynomials $P_{n}(z)$ (of first kind) and $\wt{P}_{n} (z)$ (of second kind).   We recall that $\wt{P}_n (z)$ are
  defined by equations \e{eq:RR} where $n\geq 1$ with the boundary conditions $\wt{P}_{0}(z)=0$, $\wt{P}_{1}(z)=a_{0}^{-1}$; clearly,  $\wt{P}_{n} (z)$ is a polynomial of degree $n-1$.
  
   In the case $|\beta_{\infty}| > 1$ we show that 
  the  operator $J_{0}$ is essentially self-adjoint if and only if 
   \begin{equation}
  \sum_{n =0}^{\infty} a_{n}^{-1} \big(|\beta_{\infty}|+\sqrt{\beta_{\infty}^2 -1} \big)^{2n}=\infty ;
 \label{eq:nc+}\end{equation} 
 otherwise the deficiency indices of $J_{0}$ are $(1,1)$. This result may be compared with the Berezanskii theorem (see, e.g., page 26
 in the book \cite{AKH}) stating  that   the Carleman condition \e{eq:Carl} is   necessary for the essential self-adjointness of   $J_{0}$ provided     $b_{n}=0$ and $a_{n-1}a_{n+1}\leq a_{n}^2$. As shows   Example~\ref{LCa} below, the last condition on $a_{n}$ cannot be omitted in this statement.

 %Of course for $\beta_{\infty}=1$, this condition reduces to \e{eq:Carl}.
 
 %Note that if the deficiency indices are $(1,1)$, then     the condition \e{eq:nc} holds.
 %Under the assumptions \e{eq:nc} and \e{eq:Gr} the
 %deficiency indices of the Jacobi operator $J$ are $(0,0)$ if $|\beta_{\infty}|> 1 $ and $(1,1)$ if $|\beta_{\infty}| < 1 $.

  The spectra of all  self-adjoint Jacobi operators $J$ are simple with $e_{0} = (1,0,0,\ldots)^{\top}$ being a generating vector. Therefore it is natural to define   the   spectral measure of $J$ by the relation $d\rho_{J}(\lambda)=d(E_{J}(\lambda)e_{0}, e_{0})$ where  $dE_{J}(\lambda)$      is the spectral family of the operator $J$.
  For all   extensions $J$ of the operator $J_{0}$, the  polynomials  $P_{n}(\lambda)$ are orthogonal and normalized  in the spaces $L^2 ({\Bbb R};d\rho_{J})$: 
      \[
\int_{-\infty}^\infty P_{n}(\lambda) P_{m}(\lambda) d\rho_{J}(\lambda) =\d_{n,m};
\]
as usual, $\d_{n,n}=1$ and $\d_{n,m}=0$ for $n\neq m$.
 
% Under assumptions \e{eq:Carl} and \e{eq:Gr} where $|\beta_{\infty}|< 1$, the operator
%   $J=\clos J_{0}$  has, typically, the absolutely continuous spectrum coinciding with the real axis. 

   We also note a link with a moment problem    
    \[
s_{n}= \int_{-\infty}^\infty \lambda^n d \rho (\lambda) 
\]
%label{eq:MM}\end{equation}
where $s_{n}$ are given and the measure $d \rho (\lambda) $ has to be found. 
A moment problem    is called determinate if its solution $d \rho (\lambda)$ is unique. In the opposite case it is called indeterminate.
Suppose that $s_{n}= ({\cal J}^n_{0} e_{0}, e_{0})$. Then
\[
 s_{n}= (J^n e_{0}, e_{0})= \int_{-\infty}^\infty \lambda^n d \rho_{J} (\lambda)
\]
%label{eq:MM1}\end{equation} 
for all self-adjoint extensions $J$ of the operator $J_{0}$.   It is known (see, e.g., Theorem~2 in \cite{Simon}) that the moment problem with the coefficients $ s_{n}= (J_{0}^n e_{0}, e_{0})$ is determinate if and only if the operator  $J_{0}$ is essentially self-adjoint.

   The comprehensive presentation of the results  described in this subsection  can be found in  the books \cite{AKH,   Schm} or the survey \cite{Simon}.
     We do not use  operator methods in this paper. So   the above information was given mainly to put our results into the right framework.  

% This circumstance is, however, inessential for us since we do not use any   operator methods.

 \section{Ansatz}

 In this   section, we calculate the remainder \e{eq:Grr} for the Ansatz $Q_{n} $ defined by formulas \e{eq:Ans1} or \e{eq:Ans2}. Then
  we make substitution \e{eq:Jost}.

 % Below conditions \e{eq:comp} and \e{eq:LR} are always assumed unless indicated otherwise. 

\subsection{Preliminaries}

Let us consider equation \e{eq:Jy}. Note that the values of $F_{N-1}$ and $F_{N }$ for some $N\in{\Bbb Z}_{+}$ determine the whole sequence $F_{n}$ satisfying the difference equation \e{eq:Jy}.

 Let $f=\{ f_{n} \}_{n=-1}^\infty$ and $g=\{g_{n} \}_{n=-1}^\infty$ be two solutions of equation \e{eq:Jy}. A direct calculation shows that their Wronskian
  \begin{equation}
\{ f,g \}: = a_{n}  (f_{n}  g_{n+1}-f_{n+1}  g_{n})
\label{eq:Wr}\end{equation}
does not depend on $n=-1,0, 1,\ldots$. In particular, for $n=-1$ and $n=0$, we have
 \begin{equation}
\{ f,g \} = 2^{-1} (f_{-1}  g_{0}-f_{0}  g_{-1}) \q {\rm and} \q \{ f,g \} = a_{0}  (f_{0}  g_{ 1}-f_{ 1}  g_{0})
\label{eq:Wr1}\end{equation}
(we put $a_{-1}= 1/2$).
Clearly, the Wronskian $\{ f,g \}=0$ if and only if the solutions $f$ and $g$ are proportional.

%Calculating the Wronskian \e{eq:Wr} for $n\to\infty$, we see that equation \e{eq:Jy} may have at most one (up to a   constant factor)  solution $f_{n}$ such that $f_{n}\to 0$ as $n\to\infty$.

It is convenient to introduce a notation
 \begin{equation}
x_{n}'= x_{n+1}  -x_{n}
\label{eq:diff}\end{equation}
for the ``derivative" of a sequence $x_{n}$. Note the Abel summation formula (``integration by parts"):
 \begin{equation}
\sum_{n=N }^ { M} x_{n}  y_{n}' = x_{M}  y_{M+1} - x_{N -1}  y_{N}  -\sum_{n=N } ^{M} x_{n-1}'  y_{n};
\label{eq:Abel}\end{equation}
here $M\geq N\geq 0$ are arbitrary, but we have to set $x_{-1}=0$ so that $x_{-1}'=x_{0}$.
 
   It follows from equation \e{eq:RR} that if $P_{n}  (z)$ are the orthogonal polynomials corresponding to  coefficients $( a_{n}, b_{n})$, then the polynomials $ (-1)^n P_{n}  (-z)$  correspond to the coefficients $(a_{n}, -b_{n})$. Therefore without loss of generality, we could have assumed that $\beta_{\infty}  \geq 0$.

%We  fix the branch of the analytic function $\sqrt{z^2 -1}$ of $z\in {\Bbb C}\setminus [-1,1]$ by the condition $\sqrt{z^2 -1}>0$ for $z>1$. Obviously, this function is continuous up to the cut along $[-1,1]$,    it equals $\pm i\sqrt{1-\lambda^2}$ for $z=\lambda\pm i0$, $\lambda\in (-1,1)$, and $\sqrt{z^2 -1}< 0$ for $z< -1$. Define the one-to-one, onto mapping
%$\z: {\Bbb C}\setminus [-1,1] \to {\Bbb D} $ (the unit disc) by formula \e{eq:ome}.  Then
%$$
% 2z= \z (z) +\z (z)^{-1}.
%$$
%  For $\lambda\in [-1,1]$, it is common to set $\lambda=\cos \theta$ with $\theta\in [0,\pi ]$. Then 
%$\z (\lambda\pm i0)=e^{\mp i \theta}$.

To emphasize the analogy between differential and difference operators, we often  use ``continuous" terminology (Volterra  integral equations, integration by parts, etc.) for sequences labelled by the discrete variable $n$. Below $C$, sometimes with indices,  and $c$ are different positive constants whose precise values are of no importance.

In all our constructions below, it suffices to consider the Jacobi  equation \e{eq:Jy} for large $n$ only. 

% In all estimates below 
% the values of   $|z|$ are   bounded. Then the values of    $\z=\z(z)$ are separated from $0$.

% \section{Growing coefficients}  The methods developed in the previous sections apply also to Jacobi operators with growing coefficients. 
 
 % \subsection{Ansatz}
 
  \subsection{Assumptions}
  
  In addition to \e{eq:nc}  and \e{eq:Gr}, we need some mild conditions on a regularity of behavior of the sequences $a_{n}$ and $b_{n}$ as $n\to\infty$.
   Let us set
   \begin{equation}
  \varkappa_{n}=\sqrt{\frac{a_{n+1}}{a_{n }}}, \q k_{n} = \frac{  \varkappa_{n-1}} {  \varkappa_{n}} =\frac{a_{n}  } {\sqrt{a_{n-1} a_{n+1}} } .
\label{eq:Gr4K}\end{equation}
With respect to $a_{n}$, we assume that
 \begin{equation}
 \{k_{n} -1 \}\in \ell^1 ({\Bbb Z}_{+} )
\label{eq:Gr6b}\end{equation} 
which implies also the following properties of the numbers $ \varkappa_{n} $.

\begin{lemma}\label{GK}
 Under assumption \e{eq:Gr6b} there exists a finite limit
 \begin{equation}
\lim_{n\to\infty}  \varkappa_{n} = : \varkappa_\infty,
\label{eq:Gr6a}\end{equation} 
and   $\varkappa_\infty\geq 1$ if
 condition  \e{eq:nc} is satisfied. Moreover,
  \begin{equation}
\{\varkappa_{n} '  \}\in  \ell^1 ({\Bbb Z}_{+}) 
\label{eq:Gs9b}\end{equation}
if both  \e{eq:nc} and \e{eq:Gr6b}  are true.
 \end{lemma}

   \begin{pf}
   By definition  \e{eq:Gr4K} of $k_{n}$, we have $\ln k_{n}= \ln\varkappa_{n-1} -\ln\varkappa_n $ whence
   \[
   \ln\varkappa_{n}=  \ln\varkappa_0 - \sum_{m=1}^n \ln k_m.
   \]
   It follows from \e{eq:Gr6b} that the series on the right converges which implies the existence of the limit \e{eq:Gr6a}.
   If $\varkappa_\infty <1$, then, by definition \e{eq:Gr4K} of $\varkappa_n$, we would  have $a_{n}\leq \gamma a_{n-1}$ for some $\gamma<1$ and  all $n\geq n_{0}$ if $n_{0}$ is sufficiently large. Thus, $a_{n}^{-1}\geq a_{n_{0}}^{-1}\gamma^{n_{0}}\gamma^{-n}$ so that the series in \e{eq:nc} diverges.
   
   Since
   \[
   \varkappa_{n}-\varkappa_{n-1}=\frac{ \varkappa_{n}^2}{ \varkappa_{n}+  \varkappa_{n-1} } (1+k_{n})(1-k_{n}),
   \]
   relation \e{eq:Gs9b} is a direct consequence of assumption \e{eq:Gr6b}.
      \end{pf}
      
      \begin{example}\label{GKy}. Both  conditions  \e{eq:nc} and \e{eq:Gr6b}
 are   satisfied for $a_{n}= \gamma n^p$ where $\gamma>0$, $p>1$ and for $a_{n}= \gamma x^{ n^q}$ where $x>1$, $q< 1$.  In these cases $\varkappa_\infty=1$. For $a_{n}= \gamma x^{ n}$,  conditions  \e{eq:nc} and \e{eq:Gr6b}
 are   satisfied but    $\varkappa_\infty=\sqrt{x}$. On the  contrary, condition  \e{eq:Gr6b} fails if $a_{n}= \gamma x^{ n^q}$  with $q>1$.
      \end{example}

   With respect to the coefficients $b_{n}$ or, more precisely, the ratios $\beta_{n}$ defined by \e{eq:Gr}, we assume that
     \begin{equation}
  \{\beta_{n}' \}\in \ell^1 ({\Bbb Z}_{+} ).
\label{eq:Gr8}\end{equation}

  \subsection{Construction}

Let us construct an Ansatz $Q_{n}$ satisfying relation \e{eq:Grr}  with the remainder $r_{n}  (z)\in\ell^1 ({\Bbb Z}_{+})$.
Put 
     \begin{equation}
q_{n}= Q_{n+1}  Q_{n}^{-1} .
\label{eq:Gr4L}\end{equation}
Then  \e{eq:Grr} can be rewritten as
\[
 r_{n} (z)  =   \big( \sqrt{a_{n-1}a_{n}}   \big)^{-1}
\big( a_{n-1}q_{n-1}^{-1}   + b_{n}-z  + a_{n}q_{n }\big).
\]
%label{eq:Grv}\end{equation}
Using notation  \e{eq:Gr4K} and
  introducing a  new variable
     \begin{equation}
\z_{n}= \varkappa_{n} q_{n},
\label{eq:Grz}\end{equation}
we see that
\begin{equation}
 r_{n} (z) = \zeta_{n-1}^{-1} -2\beta_{n} +
 k_{n}\zeta_{n} -2 z \alpha_{n} 
\label{eq:Grr1}\end{equation}
where
\[
\alpha_{n} =\frac{1}{2\sqrt{a_{n-1}a_{n}} } , \q  \beta_{n} =-\frac{b_{n}}{2\sqrt{a_{n-1}a_{n}} } .
\]
%label{eq:aabb}\end{equation}
Putting together relations \e{eq:Gr4L} and \e{eq:Grz} and setting $Q_{0}= a_{0}^{-1/2}$,
 we find that
 \begin{equation}
Q_{n}  = \frac{1}{\sqrt{a_{n}} } \z_{0}  \cdots \z_{n-1}.
\label{eq:ANS}\end{equation}

In principle, the numbers $\z_{n}$ can be successively determined from the equations $ r_{n} (z)=0$
which leads of course to very complicated expressions. Fortunately the construction of  $\z_{n}$ becomes quite explicit if one neglects in \e{eq:Grr1}  the terms from $\ell^1 ({\Bbb Z}_{+})$. Obviously, the term $2 z \alpha_{n} $ in the right-hand side of \e{eq:Grr1} can  be omitted if condition \e{eq:nc}  is satisfied. Furthermore,
under the assumptions below $\zeta_{n-1}^{-1} $ can be replaced by  $\zeta_{n}^{-1} $ and $k_{n}$ -- by $1 $ which allows us to define $\z_{n} $   from the equation 
 \begin{equation}
 \z_{n }+  \z_{n}^{-1}=  2 \beta_{n}   .
\label{eq:Gr5}\end{equation}
We are interested in solutions of this equation such that $| \z_{n }|\leq 1$.

%  In this paper, we consider the   case  
%  \[  \sum_{n =0}^{\infty} \alpha_n <\infty. \] 
 
Solutions of \e{eq:Gr5} are obviously given by the equalities
 \begin{equation}
\z_{n}=    \beta_n- i\sqrt{ 1- \beta_n^2}= (  \beta_n + i\sqrt{1-  \beta_n^2})^{-1} 
\label{eq:Gs3}\end{equation}
 if $| \beta_n| \leq 1 $ and
 \begin{equation}
  \z_{n}=  \sgn\beta_{n} (|  \beta_n|-  \sqrt{   \beta_n^2-1})= \sgn\beta_{n} ( | \beta_n | +   \sqrt{  \beta_n^2-1})^{-1}  
\label{eq:Gs3+}\end{equation}
 if $| \beta_n| \geq 1 $. It is convenient to set
 \begin{equation}
\theta_{n}=\arccos \beta_{n}\in [0,\pi]
\label{eq:Gk}\end{equation}
for $| \beta_n| \leq 1 $ and
 \begin{equation}
\vartheta_{n}=\arccosh |\beta_{n}|=\ln\big(|\beta_{n}|+\sqrt{\beta_{n}^2-1}\big)>0
\label{eq:Gk1}\end{equation}
for $| \beta_n| > 1 $.  Then 
 \begin{equation}
\z_{n}= e^{-i\theta_{n}} \; \mbox{for}\; | \beta_n| \leq 1\q \mbox{and}\q   \z_{n}= \sgn\beta_{n} e^{-\vartheta_{n}} \; \mbox{for}\;  | \beta_n| \geq 1.
\label{eq:Gkx}\end{equation}

It follows from 
condition \e{eq:Gr}   that
   \begin{equation}
\z_n\to \beta_{\infty} - i\sqrt{1-\beta^2_{\infty}}=:\z_\infty =e^{- i\theta_\infty}\q \mbox{if}\q | \beta_\infty| < 1 
\label{eq:zz}\end{equation}
and
   \begin{equation}
 \z_n  \to \sgn\beta_\infty (|\beta_{\infty} | -   \sqrt{\beta^2_{\infty}-1})
 =:\z_\infty = \sgn\beta_\infty e^{- \vartheta_\infty}\q \mbox{if}\q | \beta_\infty| > 1 
  \label{eq:zz+}\end{equation}
as $n\to\infty$. Here  $\theta_\infty$ and $\vartheta_\infty$ are defined by formulas \e{eq:Gk} and \e{eq:Gk1} where $n=\infty$, that is,
   \begin{equation}
\theta_\infty=\arccos \beta_\infty\in (0,\pi) \q \mbox{and}\q \vartheta_\infty= \ln\big(|\beta_\infty|+\sqrt{\beta_\infty^2-1}\big)>0 .
  \label{eq:xz+}\end{equation}
 Obviously,  $\z_\infty \in {\Bbb T}$ and $\z_\infty^2\neq 1$ if $ | \beta_\infty| < 1$ and $\z_\infty \in (-1,1)$ but $\z_\infty\neq 0$ if $ | \beta_\infty| > 1$.
For $| \beta_n| \geq 1 $ in the case $| \beta_\infty| <1 $ and for $| \beta_n| \leq 1 $ in the case $| \beta_\infty| >1 $, the numbers $\theta_{n}$ and $\vartheta_{n} $
  can be chosen in an arbitrary way; for definiteness, we set   $\theta_{n}=0$ and  $\vartheta_{n}=0$.

Now the  Ansatz \e{eq:ANS} can be written as
 \begin{equation}
Q_{n}  = \frac{1}{\sqrt{a_{n}} } e^{-i\phi_{n}}  
 \q\mbox{if}\q |\beta_{\infty}|  <1 
\label{eq:Gr1}\end{equation}
and 
 \begin{equation}
Q_{n}  = \frac{1}{\sqrt{a_{n}} } (\sgn \beta_{\infty})^n e^{- \varphi_{n}}  
 \q\mbox{if}\q |\beta_{\infty}|  >1 .
\label{eq:Gr1+}\end{equation}
Here
 \begin{equation}
 \phi_{n}=\sum_{m=0}^{n-1}\theta_{m}, \q
\varphi_{n}=\sum_{m=0}^{n-1}\vartheta_{m}
\label{eq:GRR}\end{equation}
and $\theta_n$, $\vartheta_n$  are defined by \e{eq:Gk}, \e{eq:Gk1}.  Obviously,  formulas \e{eq:Gr1}, \e{eq:Gr1+} coincide with
\e{eq:Ans1}, \e{eq:Ans2}, respectively.  

Since $\theta_{n}\to \theta_\infty$ and 
$\vartheta_{n}\to \vartheta_\infty$ as $n\to\infty$,
it follows from \e{eq:GRR} that
  \begin{equation}
\phi_{n}= \theta_\infty n + o(n) \q\mbox{and}  \q \varphi_{n}= \vartheta_\infty n + o(n) \q\mbox{as}  \q n\to\infty
\label{eq:zz1}\end{equation}
with $\theta_\infty $ and $\vartheta_\infty $  given by formulas \e{eq:xz+}.
We emphasize that the sequences $Q_{n}$ do not depend on the spectral parameter $z$. The notations $\z_{n}$ defined by \e{eq:Gs3} or \e{eq:Gs3+}  and $\z_\infty$ defined by \e{eq:zz} or \e{eq:zz+} will  often be used below.

Putting together relations \e{eq:Grr1} and \e{eq:Gr5}, we can state an intermediary result.

\begin{lemma}\label{Rr}
  Let the
     sequence $Q_{n} $ be defined by formulas \e{eq:Gr1} or \e{eq:Gr1+}.
   Then   the remainder \e{eq:Grr}  admits the  representation
\begin{equation}
 r_{n} (z)  =   \big(  \z_{n-1}^{-1}-\z_{n}^{-1}\big) + (k_{n}-1)\z_{n}-2 z \alpha_{n}
\label{eq:Gr6}\end{equation} 
where the numbers $\z_{n}$ are defined by formulas \e{eq:Gs3}  or \e{eq:Gs3+} and  $k_{n}$ is given by \e{eq:Gr4K}.
\end{lemma}

 \subsection{Estimates of the remainder}

Our goal here is to     estimate   expression \e{eq:Gr6}. 
 Under assumption \e{eq:Gr}  we have
 \begin{equation}
|1- \beta_n^2| \geq c  >0 
\label{eq:aabc}\end{equation}
  for sufficiently large $n$. It follows from definitions \e{eq:Gs3} and \e{eq:Gs3+} that
 \[
\z_{n-1}^{-1} - \z_{n}^{-1}= ( \beta_{n-1} -  \beta_{n})\Big(1- i  \frac{ \beta_{n-1}+  \beta_{n}}{\sqrt{ 1-\beta_{n-1}^2}+\sqrt{1- \beta_{n}^2}}\Big), \q |\beta_{\infty}| <1,
\]
%label{eq:Gs}\end{equation}
and
 \[
\z_{n-1}^{-1} - \z_{n}^{-1}= ( \beta_{n-1} -  \beta_{n})\Big(1 + \sgn\beta_{\infty}  \frac{ \beta_{n-1}+  \beta_{n}}{\sqrt{ \beta_{n-1}^2-1}+\sqrt{ \beta_{n}^2-1}}\Big), \q |\beta_{\infty}| >1.
\]
%label{eq:Gs+}\end{equation}
According to \e{eq:aabc} each of these identities yields  an estimate
 \begin{equation}
|\z_{n-1}^{-1} - \z_{n}^{-1}| \leq C | \beta_{n-1} -  \beta_{n}|.
\label{eq:Gs2}\end{equation}
Therefore \e{eq:Gr6} implies that
\begin{equation}
| r_{n} (z) |  \leq C    | \beta_{n-1} -  \beta_{n}|+ | k_{n} -  1| + 2 \alpha_{n}|z|
\label{eq:Gr9f}\end{equation}
where the  constant $C$ does not depend on $z$. This
leads to the following assertion.

   \begin{lemma}\label{Gr}
Let conditions  \e{eq:nc} and  \e{eq:Gr}   be satisfied. Then   estimate \e{eq:Gr9f} for the remainder \e{eq:Gr6} is true for all $z\in{\Bbb C}$. Under additional
 assumptions \e{eq:Gr6b}  and \e{eq:Gr8}, 
we have
\[
\{ r_{n} (z) \} \in \ell^1 ({\Bbb Z}_{+}) .
\]
%label{eq:Gr9}\end{equation}
\end{lemma}
     
%Condition \e{eq:Gr6b} implies some restrictions on the behavior as $n\to\infty$ of the numbers $\varkappa_{n}$ defined in \e{eq:Gr4K}.

%It is convenient to make a weak additional assumption 
%\begin{equation}
%\lim_{n\to\infty}  \frac{a_{n+1}} {a_{n}} =1
%\label{eq:Gr6a}\end{equation} 
%which excludes a super-power growth of $a_{n}$ as $n\to\infty$.
% Then   $k_{n}\to 1$ as $n\to\infty$,  but we need a stronger condition  \e{eq:Gr6b}. Note that it is satisfied if
% $a_{n}=n^p$ for some $p>0$.

%[ R_{n} (z)  =  ( \frac{1} {\z_{n-1}} - \frac{1} {\z_{n}} )+(\frac{\sqrt{a_{n-1} a_{n+1}} }{a_{n}  }-1) \z_{n-1}^{-1}.\]
%Of course if $k_{n} =1$, then $\varkappa_{n} =1$, expression \e{eq:Gs3} reduces to \e{eq:aa1} and
%\begin{equation} R_{n} (z)  =  \frac{a_{n-1}} {a_{n}}  \frac{1} {\z_{n-1}} - \frac{1} {\z_{n}}. \label{eq:Gr6old}\end{equation}
% However this choice does not work in the case \e{eq:Gr}. Indeed, if, for example,   $a_{n}=n^p$ for some $p>0$, then $a_{n+1}a_{n}^{-1}=1+ p n^{-1}+ O(n^{-2})$ and the remainder \e{eq:Gr6old} does not belong to $\ell^1 ({\Bbb Z}_{+})$.

\subsection{Multiplicative substitution}

Let the sequence
$ Q_{n} $ be given by formulas \e{eq:Gr1} or \e{eq:Gr1+}.
We are looking for solutions 
  $f_{n} (z)$  of the difference equation \e{eq:Jy} satisfying the condition 
   \begin{equation}
f_{n} (z) = Q_{n} (1+ o(1)), \q n\to\infty.
\label{eq:Gs1}\end{equation}
The uniqueness of such solutions is almost obvious.

\begin{lemma}\label{uniq}
Equation \e{eq:Jy} may have only one solution   $f_{n} (z)$ satisfying condition \e{eq:Gs1}.
 \end{lemma}

\begin{pf}
Let $\mathtt{f}_{n} (z)$ be another solution of \e{eq:Jy} satisfying  \e{eq:Gs1}. It follows from \e{eq:Gr1} or \e{eq:Gr1+} that the Wronskian \e{eq:Wr} of these solutions
calculated for $n\to\infty$ equals
\[
\{f, \mathtt{f} \}=a_{n} Q_{n} Q_{n+1} o(1)= o(1)\varkappa_{n}^{-1} \begin{cases}e^{-i(\phi_{n-1}+\phi_{n})},\q |\beta_{\infty}| <1\\
\sgn\beta_{\infty}e^{-\varphi_{n-1}-\varphi_{n}},\q |\beta_{\infty}| >1\end{cases} 
\]
so that $\{f, \mathtt{f}\}=0$ according to \e{eq:Gr6a}. Thus $\mathtt{f} =Cf$ where $C=1$ by virtue again of condition \e{eq:Gs1}.
\end{pf}

\begin{remark}\label{uniqq}
For the calculation of the Wronskian $\{f,  \mathtt{f} \}$, it is essential that the power of $a_{n}$ in \e{eq:ANS} equals $-1/2$.
 \end{remark}

For construction of $f_{n}  (z)$, we will reformulate the problem    introducing a  sequence
\begin{equation}
 u_{n} (z)=  Q_{n} ^{-1}  f_{n} (z), \q n\in {\Bbb Z}_{+}.
\label{eq:Gs4}\end{equation}
 Then \e{eq:Gs1} is equivalent to the condition  
\[
\lim_{n\to\infty} u_{n} (z)=   1.
\]
%label{eq:A12a}\end{equation}

Next, we derive a difference equation for $ u_{n} (z)$.

\begin{lemma}\label{Gs}
Let $z\in\Bbb C$, let $\z_{n}$ be defined by formulas \e{eq:Gs3} or \e{eq:Gs3+}, and let the remainder $r_{n} (z) $ be given by formula \e{eq:Gr6}. Then
 equation  \e{eq:Jy} for a sequence $ f_{n} (z)$ is equivalent to the equation
\begin{equation}
 k_{n} \z_{n} ( u_{n+1} (z)- u_{n} (z)) -    \z_{n-1}^{-1}   ( u_{n} (z)- u_{n-1} (z))=-   r_{n} (z) u_{n} (z), \q n\in {\Bbb Z}_{+},
\label{eq:Gs5}\end{equation}
for  sequence  \e{eq:Gs4}. 
 \end{lemma}

\begin{pf}
Substituting  expression $f_{n}  = Q_{n}  u_{n} $ into  \e{eq:Jy} and using the  equality 
\[\frac{Q_{n+1}}{Q_{n }}= \sqrt{\frac{a_{n}}{a_{n +1}}} \z_{n}, 
\]   we see that
\begin{multline*}
( \sqrt{a_{n-1}a_{n}}Q_{n} )^{-1}\Big(  a_{n-1} f_{n-1} + ( b_{n} -z)f_{n} + a_{n} f_{n+1}\Big)
\\=
\sqrt{\frac{a_{n-1} }{a_{n}}}\frac{Q_{n-1} }{Q_{n}} u_{n-1}- 2 (\alpha _{n} z  + \beta_{n}) u_{n}+
\sqrt{\frac{a_{n} }{a_{n-1}}}  \frac{Q_{n+1} }{Q_{n}} u_{n+1}
\\  
 =    \z_{n-1}^{-1} u_{n-1} -2   (\alpha _{n} z  + \beta_{n})  u_{n} + k_{n}   \z_{n} u_{n+1}   .
% \label{eq:Gs6}
  \end{multline*}
 In view of  \e{eq:Gr5} the right-hand side here equals
  \[  
   k_{n}   \z_{n}( u_{n+1} -u_{n}) -\z_{n-1}^{-1} (u_{n} -u_{n-1} )+ r_{n} u_{n} 
 \]
 with $r_{n}$ defined by \e{eq:Gr6}.
 Therefore    the equations \e{eq:Jy}  and \e{eq:Gs5} are equivalent.
 \end{pf}
 
  \section{Modified Jost solutions }

     In this  section,
  we reduce
 the Jacobi difference equation  \e{eq:Jy} for Jost solutions $f_{n}  (z)$   to a Volterra equation for functions $u_{n}  (z)$ defined by formula  \e{eq:Gs4} and satisfying the condition $u_{n}  (z)\to 1$ as $n\to\infty$.  
 Solutions $u_{n}  (z)$ of the Volterra  equation  can be constructed by iterations. At this point there  is almost no difference between  the cases $|\beta_{\infty}| <1$ and $|\beta_{\infty}| >1$. Finally, formula  \e{eq:Gs4}  yields asymptotics of $f_{n}  (z) $ as $n\to\infty$.

\subsection{Volterra integral equation}

The sequence $u_{n}  (z)$ satisfying the difference equation \e{eq:Gs5} will be constructed as a solution of an appropriate ``Volterra integral" equation.  We set
 \begin{equation}
\sigma_{n}=    \z_{n} \z_{n-1}   ,\q   S_{n}=   \sigma_{0} \sigma_{1}  \cdots \sigma_{n-1},\q n\geq 1,
\label{eq:Gs7}\end{equation}
and, for $ n<m $,
\begin{equation}
G_{n,m}=  -  (\varkappa_{m-1} \z_{m})^{-1} S_{m+1}  \sum_{p=n+1}^m \varkappa _{p-1}  S_{p}^{-1} = -  (\varkappa_{m-1} \z_{m})^{-1}\sum_{p=n+1}^m \varkappa_{p-1}  \sigma_{p}  \cdots \sigma_{m}   .
\label{eq:Gp6}\end{equation}
Note  that the kernel $G_{n,m}$ does not depend on $z$.
Let us consider an equation
\begin{equation}
u_{n} (z)= 1+ \sum_{m=n+1}^\infty G_{n,m} r_m (z) u_{m} (z)
\label{eq:Gp5}\end{equation}
where the sequence $r_m (z) $ is defined by formula \e{eq:Gr6}.

%We suppose that $z$ belongs to a compact subset of $\Bbb C$, $N$  is sufficiently large  and $n\geq n_{0}$. 

Our first goal is to estimate the ``kernel" $G_{n,m}$. This is quite straightforward in the case $|\beta_{\infty}|>1$.

 \begin{lemma}\label{GS2+}
Let  the assumptions  \e{eq:nc}, \e{eq:Gr} where $|\beta_{\infty}|>1$ and  \e{eq:Gr6b}  be satisfied.   Then
 kernel \e{eq:Gp6} is bounded uniformly in $m>n\geq 0$:
\begin{equation}
|G_{n,m}|\leq    C<\infty.
\label{eq:AbG}\end{equation}
 \end{lemma}

 \begin{pf}
  By virtue of \e{eq:zz+}, we have 
 \begin{equation}
|\sigma_n|\leq s_{\infty} <1
\label{eq:ss+}\end{equation}
  if $\z_{\infty}^2 < s_{\infty}$ and $n$ is sufficiently large, say, $n\geq n_{0}$. Moreover, $\{\varkappa_{n}\}\in \ell^2 ({\Bbb Z}_{+})$ and $\{\varkappa_{n}^{-1}\}\in \ell^2 ({\Bbb Z}_{+})$ according to Lemma~\ref{GK}.
    Therefore it follows from  \e{eq:Gp6} that
\[
|G_{n,m}| \leq  C  
  \sum_{p=n+1}^m s_{\infty}^{m-p+1} =C  \frac{1- s_{\infty}^{m-n}}  
  {s_{\infty}^{-1}  -1 }   \leq C  \frac{ s_{\infty} }  
  {1-s_{\infty}  }   ,\q n<m  ,
\]
% label{eq:Gp6+}\end{equation}
if $n\geq n_{0}$. 
 
  Let now $n< n_{0}$. If $n< m\leq n_{0}$, then the sum \e{eq:Gp6}  consists of at most $n_{0}$ terms. If $m>n$, then 
\[
G_{n,m}=G_{n,n_{0}}+G_{n_{0},m}
\]
is a sum of two bounded terms.
   \end{pf}

In the case $|\beta_{\infty}|<1$ we have to ``integrate by parts".

 \begin{lemma}\label{GS2}
Let  the assumptions  \e{eq:nc}, \e{eq:Gr} where $|\beta_{\infty}|<1$ and  \e{eq:Gr6b}, \e{eq:Gr8}  be satisfied.  Then
 kernel \e{eq:Gp6} is bounded uniformly in $m>n\geq 0$, that is, estimate
\e{eq:AbG} holds.
 \end{lemma}

 \begin{pf}
 As in the previous lemma, we can suppose   that $n\geq n_{0}$ where  $n_{0}$  is sufficiently large.
 Since $|\z_{p}|= 1$ according to \e{eq:Gs3}, it follows from 
  definition \e{eq:Gs7} that
  \begin{equation}
|S_{m+1} S_{p}^{-1}| =|\sigma_{p} \cdots \sigma_{m} | =  1,\q p \leq m.
\label{eq:ssc}\end{equation}
 Relation \e{eq:zz} where $\z_{\infty}^2\neq 1$ implies that
 \begin{equation}
|\sigma_n-1|\geq c>0 .
\label{eq:ss}\end{equation}
By definitions \e{eq:diff} and \e{eq:Gs7} we have
\[
(S_{p } ^{-1})'=S_{p+1} ^{-1}  - S_{p } ^{-1}= ( \sigma_{p}^{-1} -1 )  S_{p}^{-1}, 
\]
and hence integrating by parts (that is, using formula \e{eq:Abel}), we find that
\begin{multline}
\sum_{p=n+1 }^ { m} \varkappa_{p-1} S_{p}^{-1} = \sum_{p=n+1 }^ { m} \varkappa_{p-1}  ( \sigma_{p}^{-1}  -1)^{-1}(S_{p}^{-1})'
 \\
=\varkappa_{m-1} ( \sigma_{m}^{-1}  -1)^{-1}S_{m+1}^{-1}-  \varkappa_{n-1} ( \sigma_{n}^{-1}  -1)^{-1}S_{n+1}^{-1}-\sum_{p=n+1 }^ { m} (\varkappa_{p-2}  (\sigma_{p-2}^{-1}  -1)^{-1})'S_{p}^{-1} .
\label{eq:AbelG}\end{multline}
Note that  
\[ 
(( \sigma_{p-1}^{-1}  -1)^{-1})'=  ( \sigma_{p-1}  -1)^{-1}( \sigma_{p}  -1)^{-1} \sigma_{p-1}' 
\]
where $\sigma_{p}'\in  \ell^1 ({\Bbb Z}_{+})$ by  virtue of \e{eq:Gs2} and \e{eq:Gr8}. Using also \e{eq:Gs9b} and  \e{eq:ss},  we see that
this sequence belongs to $\ell^1 ({\Bbb Z}_{+})$ and hence
\[
(\varkappa_{p-1}  (\sigma_{p-1}^{-1}  -1)^{-1})'\in  \ell^1 ({\Bbb Z}_{+}).
\]
Let us multiply identity \e{eq:AbelG} by $S_{m+1} $.
According to \e{eq:ssc} and \e{eq:ss} all three terms in the right-hand side of the equality obtained are bounded   for $n\geq n_{0}$. 
   \end{pf}

    \subsection{Successive approximations}
    
    Let us come back to the Volterra equation \e{eq:Gp5}. Lemma~\ref{Gr}  gives a necessary estimate on the sequence 
    $r_n (z) $ and   Lemmas~\ref{GS2+} or \ref{GS2} show that the kernels $G_{n,m}$ are bounded. This
    allow us to estimate iterations of equation \e{eq:Gp5} and then solve it.

%    standardly  solve the Volterra equation \e{eq:Gp5} by iterations. 

  \begin{lemma}\label{GS3p}
  Let  the assumptions  \e{eq:nc}, \e{eq:Gr} as well as  \e{eq:Gr6b},  \e{eq:Gr8}  be satisfied. 
  % We also  suppose that  the conditions of Lemma~\ref{GS2+} are satsified in the case $|\beta_{\infty}| >1$ and that the conditions of Lemma~\ref{GS2} are satisfied in the case $|\beta_{\infty}| <1$.  
  Set $u^{(0)}_n =1$   and 
  \begin{equation}
 u^{(k+1)}_{n}(z)= \sum_{m=n +1}^\infty G_{n,m} r_{m} (z) u^{(k )}_m (z),\q k\geq 0,
\label{eq:W5}\end{equation}
for all $n\in {\Bbb Z}_{+}$. Then the estimates
 \begin{equation}
| u^{(k )}_{n} (z) |\leq \frac{C^k}{k!} \big(\sum_{m=n+1}^\infty |r_{m} (z)|\big)^k,\q \forall k\in{\Bbb Z}_{+}.
\label{eq:W6s}\end{equation}
are true for all sufficiently large $n$ with the same constant $C $ as in Lemmas~\ref{GS2+} or \ref{GS2}. 
\end{lemma}

 \begin{pf}
  Suppose that \e{eq:W6s} is satisfied for some $k\in{\Bbb Z}_{+}$. We have to check 
 the same estimate (with $k$ replaced by $k+1$ in the right-hand side)  for $ u^{(k+1)}_{n}$.  
 Set
 \[
 { R}_{m}= \sum_{p=m +1}^\infty  | r_p | .
 \]
  According to definition \e{eq:W5}, it  follows from estimate \e{eq:AbG} and  \e{eq:W6s} that
   \begin{equation}
| u^{(k +1)}_{n} |\leq \frac{C^{k+1}}{k!}  \sum_{m=n +1}^\infty  | r_m | { R}_{m}^k.
\label{eq:V7}\end{equation}
Observe that
 \[
{ R}_{m}^{k+1}+ (k+1)   | r_{m}|  { R}_{m}^k  
\leq
{ R}_{m-1}^{k+1},
\]
and hence, for all $N\in{\Bbb Z}_{+}$,
   \[
 (k+1)  \sum_{m=n +1}^N  | r_m | {  R}_{m}^k 
 \leq 
 \sum_{m=n +1}^N  ( {  R}_{m-1}^{k+1}-{  R}_{m}^{k+1})
= { R}_{n}^{k+1}-{  R}_{N}^{k+1}\leq  {  R}_{n}^{k+1}.
 \]
Substituting this bound into   \e{eq:V7}, we obtain estimate \e{eq:W6s} for $u^{(k +1)}_{n}$.
    \end{pf}

Now we are in a position to solve equation \e{eq:Gp5} by iterations.

% Now we are in a position to consider equation \e{eq:Gp5}. For the proof of the next result, we combine Lemmas~\ref{Gr}, \ref{GS2+} and \ref{GS2}.
    
    %Let us now consider equation \e{eq:Gp8}. First, we have to get rid of the assumption \e{eq:Gp4}.

  \begin{theorem}\label{GS3}
  Let the assumptions  \e{eq:nc},  \e{eq:Gr}   as well as  \e{eq:Gr6b},  \e{eq:Gr8} be satisfied.
  %  In the case $|\beta_{\infty}| <1$ we impose additionally   condition  \e{eq:Gs9b}. 
  Then equation \e{eq:Gp5} has a unique bounded solution $u_{n} (z)$. Moreover,
  \begin{equation}
 |u _{n} (z) -1 |\leq C \varepsilon_{n},\q n\geq 0,
\label{eq:Gp9}\end{equation}
where
 \begin{equation}
\varepsilon_{n}= \sum_{m=n+1}^\infty ( | \beta_{m-1} -  \beta_{m}|+ | k_{m} -  1| +  \alpha_{m}|z|)
\label{eq:Gp10}\end{equation} 
and the constant $C$   does not depend on $z\in {\Bbb C}$.
For all $n\in{\Bbb Z}_{+}$, the functions $u_{n}(z)$    are  entire functions of $z\in \Bbb C$.
\end{theorem}

 \begin{pf}  Set
     \begin{equation}
   u_{n} =\sum_{k=0}^\infty u^{(k)}_{n} 
\label{eq:W8}\end{equation}
where $u^{(k)}_{n}$ are defined by recurrence relations \e{eq:W5}.
Estimate \e{eq:W6s} shows that this series is absolutely convergent. Using the Fubini theorem to interchange the order of summations in $m$ and $k$, we see that
\[
   \sum_{m=n+1}^\infty G_{n,m}     r_{m}  u_{m}    =   \sum_{k=0}^\infty\sum_{m=n+1}^\infty G_{n,m}   r_{m}  u_{m}^{(k)} =- \sum_{k=0}^\infty  u_{n}^{(k+1)}=1- \sum_{k=0}^\infty  u_{n}^{(k)}=1- u_{n}.
\]
This is equation  \e{eq:Gp5} for sequence \e{eq:W8}. Estimate \e{eq:Gp9}  also follows from \e{eq:W6s}, \e{eq:W8}. 

According to \e{eq:Gr6}  the remainder $r_{m}(z)$ and hence
 the kernels $G_{n,m} r_{m} (z)$ are linear functions of $z$. Therefore recurrence arguments show that all 
  successive approximations $u_{n}^{(k)} (z)$     are  analytic functions of $z\in \Bbb C$. The same assertion is of course true for the series   \e{eq:W8}.
\end{pf} 

  % Under the assumptions of Lemmas~\ref{GS2+} or \ref{GS2}
 
It turns out that the construction above yields a solution of the difference equation  \e{eq:Gs5}.
 
  \begin{lemma}\label{GS4}
  Let $r_{n}(z)$ and $G_{n,m}$ be given by formulas  \e{eq:Gr6} and  \e{eq:Gp6}, respectively.
Then  a solution $u_{n}  (z)$ of   integral  equation   \e{eq:Gp5}   satisfies also the difference equation  \e{eq:Gs5}.
 \end{lemma}
 
  \begin{pf}
  It follows from \e{eq:Gp5}    that
  \begin{equation}
 u_{n+1} - u_{n}=   \sum_{m=n+2}^\infty (G_{n+1, m}-G_{n, m})r_m u_{m} -
 G_{n, n+1} r_{n+1} u_{n+1}.
  \label{eq:A17Ma}\end{equation}
  Since according to  \e{eq:Gp6}  
\[
G_{n+1, m}-G_{n, m}=- \varkappa_{n} (\varkappa_{m-1} \z_{m} )^{-1} S_{n+1}^{-1} S_{m+1} \q \mbox{and}\q
G_{n, n+1}=  \z_{n+1}^{-1}  S_{n+2} S_{n+1}^{-1} ,
\]
equality \e{eq:A17Ma} can be rewritten as
\begin{equation}
 \varkappa_{n}^{-1} (u_{n+1} - u_{n} )=  - \sum_{m=n+1}^\infty      S_{n+1}^{-1} S_{m+1} (\varkappa_{m-1} \z_{m} )^{-1}  r_{m} u_{m} .
  \label{eq:A17Mb}\end{equation}
Putting together this equality with the same equality for $n+1$ replaced by $n$, we see that
 \begin{multline*}
  \varkappa_{n}^{-1} ( u_{n+1} - u_{n}) - \sigma_{n}^{-1}  \varkappa_{n-1}^{-1} ( u_{n} - u_{n-1})\\=  \sum_{m=n+1}^\infty    S_{n+1}^{-1} S_{m+1}  r_{m} (\varkappa_{m-1} \z_{m} )^{-1} u_{m}-
\sigma_{n}^{-1} \sum_{m=n}^\infty  S_{n}^{-1} S_{m+1}    r_{m} (\varkappa_{m-1} \z_{m} )^{-1} u_{m}.
  \end{multline*}
Since  $S_{n+1}=\sigma_{n} S_{n}$, the right-hand side here equals $ -(\varkappa_{n-1} \z_n )^{-1}r_{n}u_{n}$, and hence the equation obtained coincides with \e{eq:Gs5}.  
   \end{pf}
   
   The above arguments show also that the functions $u_{n} (z)$ are of minimal exponential type.

  \begin{lemma}\label{Gexp}
  Under the assumptions of Theorem~\ref{GS3},  for an arbitrary $\varepsilon>0$ and some constants $C_{n}(\varepsilon) $ $($that   do not depend on $z\in{\Bbb C})$, every function $u_{n} (z)$ satisfies an estimate
  \begin{equation}
 |u _{n}(z)  |\leq C_{n}(\varepsilon) e^{\varepsilon |z|}.
\label{eq:ss1}\end{equation}
\end{lemma}

  \begin{pf}
According to   \e{eq:Gp9}  and \e{eq:Gp10} we have an estimate
   \begin{equation}
 |u _m(z)  |\leq 1+ \varepsilon_m(1+|z|)\leq  2 e^{\varepsilon_m |z|}
\label{eq:ss2}\end{equation}
where $\varepsilon_m\to 0$ as $m\to\infty$.
On the other hand, it follows from equation  \e{eq:Jy} for function  \e{eq:Jost} that 
   \begin{equation}
|u_{n} (z)| \leq C_{n} (1+ |z|) (|u_{n+1} (z)|+|u_{n+2} (z)|)\leq\cdots
 \leq C_{n,k} (1+ |z|)^k (|u_{n+k} (z)|+|u_{n+k+1} (z)|)
\label{eq:ss2+}\end{equation}
for every $k=1,2, \ldots$. For a given $\varepsilon>0$, choose $k$ such that $2\varepsilon_{n+k }\leq \varepsilon$, 
$2\varepsilon_{n+k+1 }\leq \varepsilon$. Then
putting  estimates \e{eq:ss2} and \e{eq:ss2+} together, we see that 
   \[
 |u _n(z)  |\leq 4 C_{n,k} (1+ |z|)^k e^{\varepsilon |z|/2} .
\]
Since $(1+ |z|)^k\leq c_{k} (\varepsilon) e^{\varepsilon |z|/2}$,
this proves \e{eq:ss1}.
    \end{pf}

In the case $|\beta_{\infty}|>1$, we   also need 
      estimates on $u_{n}'$.
      
       \begin{lemma}\label{GS'}
       Let $|\beta_{\infty}|>1$.  Under the assumptions of Theorem~\ref{GS3}, we have  
 $\{u_{n}' \} \in  \ell^1 ({\Bbb Z}_{+}) $.
\end{lemma}

  \begin{pf} 
  Let us proceed from expression   \e{eq:A17Mb} for $u_{n+1} - u_{n}$.
  It follows from \e{eq:Gs7} and \e{eq:ss+}  that    
   \[
  |S_{n+1}^{-1} S_{m+1}|= | \sigma_{n+1}\cdots \sigma_{m}| \leq s_{\infty}^{ m-n}
  \]
  for sufficiently large $n$ and $m\geq n$ whence
  \[
  | u_{n+1} -u_{n}| \leq C \sum_{m=n+1}^\infty s_{\infty}^{m-n} | r_{m} u_{m}|.
  \]
  So we only have to take the sum over $n$ here. Then we use that $s_{\infty}  <1$ 
  and  $\{r_m u_m \} \in  \ell^1 ({\Bbb Z}_{+}) $.
\end{pf}

     %   \subsection{Asymptotic behavior}
     
          \subsection{Jost solutions}
  
  Now we are in a position to construct solutions of  the difference equation \e{eq:Jy} with asymptotics \e{eq:Gs1} as $n\to\infty$. Recall    that the sequence $Q_{n}$ is defined by formulas \e{eq:Gr1} or \e{eq:Gr1+}. According to Lemma~\ref{uniq},  equation \e{eq:Jy} may have only one solution with such asymptotics.
  By Lemma~\ref{GS3},
      equation \e{eq:Gp5} has a  solution $u_{n}(z)$ satisfying estimate
\e{eq:Gp9}  with the remainder $\varepsilon_{n}$  given by formula  \e{eq:Gp10}. By Lemma~\ref{GS4}    the difference equation \e{eq:Gs5} is also true for this solution.
 It follows from  Lemma~\ref{Gs} that
   \begin{equation}
   f_{n}( z )= Q_{n}  u_{n} (z )
   \label{eq:JostG}\end{equation}
 satisfies equation \e{eq:Jy}. 
   
   This leads to the following results. We state them separately for the cases $|\beta_{\infty}| < 1$ and $|\beta_{\infty}| > 1$.

    \begin{theorem}\label{GSS}
      Let the assumptions  \e{eq:nc},  \e{eq:Gr}  with $|\beta_{\infty}|< 1$ as well as  \e{eq:Gr6b},  \e{eq:Gr8} be satisfied.
   Then    the equation \e{eq:Jy} has a solution $\{f_{n}( z )\}$ with asymptotics
    \begin{equation}
f_{n}(  z )  = \frac{1}{\sqrt{a_{n}} } e^{-i\phi_{n}} \big(1 + O( \varepsilon_{n})\big) , \q n\to \infty,
\label{eq:A22G}\end{equation}
where $ \phi_{n}$ is given by formula \e{eq:Grf}. 
Asymptotics \e{eq:A22G} is  uniform in $z$ from compact subsets    of the complex plane $  \Bbb C $.
For all $n\in {\Bbb Z}_{+}$, the functions $f_{n}( z )$ are entire functions of $z\in  \Bbb C$ of minimal exponential type.
 \end{theorem}

   By analogy with the continuous case,   the sequence $\{f_{n}( z )\}_{n=-1}^\infty$   will be called
  the  (modified) Jost solution of equation \e{eq:Jy}. Additionally, 
we define the conjugate Jost solution by the formula
   \begin{equation}
\tilde{f}_{n}( z)=\ov{f_n( \bar{z })}.
 \label{eq:AcGG}\end{equation}
 It also satisfies equation \e{eq:Jy} because the coefficients $a_{n}$ and $b_{n}$ are real, and it has the asymptotics
   \begin{equation}
\tilde{f}_{n}(  z )= \frac{1}{\sqrt{a_{n}} }  e^{ i\phi_{n}} \big(1 + O( \varepsilon_{n})\big), \q n\to\infty,
\label{eq:A2CG}\end{equation}

Using  asymptotic formulas \e{eq:A22G}, \e{eq:A2CG} and  definition  \e{eq:Gk}, we can calculate the Wronskian of the solutions $f(z)$, $\tilde{f} (  z )$: 
  \[
\{ f(z), \tilde{f} (  z )\} =\lim_{n\to\infty}  \varkappa_{n}^{-1} (e^{i\theta_{n}} - e^{-i\theta_{n}})=\varkappa_{\infty}^{-1} (e^{i\theta_\infty} - e^{-i\theta_\infty}).
\]
Since according to \e{eq:Gs3}, \e{eq:Gkx}
\[
e^{i\theta_\infty}=\beta_{\infty}+ i \sqrt {1-\beta_{\infty}^2},
\]
it follows that
  \begin{equation}
\{ f(z), \tilde{f} (  z )\} =  2i \varkappa_{\infty}^{-1} \sqrt {1-\beta_{\infty}^2}\neq 0.
\label{eq:A2CG1}\end{equation}
  Thus, in  the case $|\beta_{\infty}| < 1$ for all  $z\in{\Bbb C}$, we have two linearly independent oscillating solutions $f_{n}(  z )$ and $\tilde{f}_{n}(  z )$. Under assumption \e{eq:nc} both of them belong to $\ell^2 ({\Bbb Z}_{+})$.  Note that the amplitude factors $a_{n} ^{-1/2}$ in \e{eq:A22G} and \e{eq:A2CG} are quite natural (cf. Remark~\ref{uniqq}) because their product should  cancel $a_{n}$ in \e{eq:Wr}.

Similarly, in the case $|\beta_{\infty}| > 1$, we have the following result.

 \begin{theorem}\label{GSS+}
      Let the assumptions  \e{eq:nc},  \e{eq:Gr}  with $|\beta_{\infty}| > 1$ as well as  \e{eq:Gr6b},  \e{eq:Gr8} be satisfied.
       Then    the equation \e{eq:Jy} has a solution $\{f_{n}( z )\}$ with asymptotics
    \begin{equation}
f_{n} (z) = \frac{1}{\sqrt{a_{n}} } (\sgn \beta_{\infty})^n e^{- \varphi_{n}} \big(1 + O( \varepsilon_{n})\big), \q n\to\infty, 
  \label{eq:A22G+}\end{equation}
where $\varphi_{n}$ is given by formula \e{eq:Grf1}.
Asymptotics \e{eq:A22G+} is  uniform in $z$ from compact subsets    of the complex plane $  \Bbb C $.
For all $n\in {\Bbb Z}_{+}$, the functions $f_{n}( z )$ are entire functions of $z\in  \Bbb C$ of minimal exponential type.
 \end{theorem}
 
 According to \e{eq:xz+} and \e{eq:zz1} 
 the solution \e{eq:A22G+} tends to zero exponentially as $n\to\infty$. It will   also be called the (modified)  Jost solution  of equation \e{eq:Jy}. However, in contrast to the case $|\beta_{\infty}| < 1$,
for $|\beta_{\infty}| > 1$ the construction of  Theorem~\ref{GSS+} yields only one  solution $f_{n}(  z )$ of the Jacobi equation \e{eq:Jy}.

%    \begin{equation}
%\varphi_{n}=n\ln\big(|\beta_\infty|+\sqrt{\beta_\infty^2-1}\big)+ o(n)\q\mbox{as}\q n\to\infty
%\label{eq:Gk2}\end{equation}
  
   \begin{example}\label{ex1}
   Suppose that (for sufficiently large $n$)
      \[
    a_{n}=\gamma n^p \q \mbox{and}\q  b_{n}=\d n^q
    \]
    where $\gamma>0$, $p>1$. 
    The assumptions of Theorem~\ref{GSS} are satisfied if  
    either $q<p$ or $q=p$ but $|\d| <2\gamma$. The assumptions of Theorem~\ref{GSS+} are satisfied if $q=p$  and $|\d|  > 2\gamma$. 
 \end{example}
  
%  \begin{example}\label{ex2}
%    The assumptions of Theorem~\ref{GSS+} are satisfied if (for sufficiently large $n$)
  %  \[   a_{n}=\gamma n^p \q \mbox{and}\q  b_{n}=\d n^p  \]
  %  where $\gamma>0$, $p>1$ and   $|\d|  > 2\gamma$. 
% \end{example}

  \begin{remark}\label{const}
The definitions \e{eq:A22G} and \e{eq:A22G+} of the Jost solutions remain unchanged (up to  insignificant constant factors) if the coefficients $b_{n}$ are replaced by $\tilde{b}_{n} =b_{n}+c$ for some $c\in{\Bbb R}$. Indeed, for the corresponding numbers \e{eq:Gr}, we have $\{\tilde{\beta}_{n}  -   \beta_{n}\} \in \ell^1 ({\Bbb Z}_{+})$ and hence the differences $\tilde{\phi}_{n} -\phi_{n}$ and $\tilde{\varphi}_{n} -\varphi_{n}$ of the phases \e{eq:Grf} and \e{eq:Grf1} have finite limits as $n\to\infty$.
 \end{remark}

 %    By analogy with the continuous case,   the sequence $\{f_{n}( z )\}_{n=-1}^\infty$   will be called
 % the  (modified) Jost solution of equation \e{eq:Jy} and $\Omega(z)$ will be called   the (modified) Jost function.  
          
 %    The following result is a direct consequence of Theorem~\ref{GSS}.

%\begin{corollary}\label{JOSTG}
%The Jost function $\Omega(z)$   depends analytically on $z\in\Bbb C$ and   satisfies estimates \e{eq:ss1}. 
%\end{corollary}

%Then the congugate Jost function $\wt\Omega(z)$ equals
%  \begin{equation}
%\wt\Omega (z):=  \{ P(z), \tilde{f}(z)\} = - 2^{-1}\tilde{f}_{-1}(z).
%\label{eq:WRGc}\end{equation}

   \section{Orthogonal polynomials}

   % The case $|\beta_{\infty}|<1$ is considered in Sect. 4.1, and the case $|\beta_{\infty}|>1$ is considered in Sect. 4.2. 
 
      \subsection{ Small diagonal elements}

Here we suppose that $|\beta_{\infty}|<1$. Recall that the Jost solution $f(z)=\{f_n (z)\}$ was constructed in Theorem~\ref{GSS} and $\tilde{f}  (z)$ was  defined by relation \e{eq:AcGG}.
By virtue of \e{eq:A2CG1}, an arbitrary solution $F(z)=\{F_{n}(z)\} $ of the Jacobi equation  \e{eq:Jy} is a linear combination of the Jost solutions $f(z)$ and $ \tilde{f} (  z )$, that is 
  \begin{equation}
F_{n}(z)= k_{+} (z)f_{n}(z) + k_{-} (z) \tilde{f}_{n}(  z ),
\label{eq:LC}\end{equation}
where the constants can be expressed via the Wronskians:
  \begin{equation}
  k_{+}(z)= \varkappa_{\infty} \frac{\{F(z), \tilde{f} (  z )\}} {2i \sqrt {1-\beta^2}}, \q
   k_{-}(z)=- \varkappa_{\infty}\frac{\{F(z), f(  z )\} }{2i \sqrt {1-\beta^2}} .
\label{eq:LC1}\end{equation}

According to \e{eq:LC}     the following result is a direct consequence of Theorem~\ref{GSS}. Recall that the phase $\phi_{n}$ is defined by \e{eq:Grf} and the remainder $\varepsilon_{n}$ is given by \e{eq:Gp10}.

 \begin{theorem}\label{LC}
    Let the assumptions of Theorem~\ref{GSS} be satisfied. Choose some $z\in {\Bbb C}$.
 Then  an arbitrary solution $F_{n}   $ of the Jacobi equation  \e{eq:Jy} has  asymptotics
  \begin{equation}
F_{n} = a_{n} ^{-1/2} \big(k_{+} e^{-i\phi_{n}} + k_{-} e^{i\phi_{n}}\big) \big(1 + O( \varepsilon_{n})\big), \q n\to\infty,
\label{eq:LC2}\end{equation} 
for some $k_{\pm }\in {\Bbb C}$. Conversely, for arbitrary $k_{\pm}\in {\Bbb C}$, there exists a solution $F_{n}$  of the  equation  \e{eq:Jy} with asymptotics \e{eq:LC2}. 
 \end{theorem}

    \begin{corollary}\label{LCx}
Since all solutions of the Jacobi equation \e{eq:Jy} are in $\ell^2 ({\Bbb Z}_{+})$,  under the assumptions of Theorem~\ref{GSS}  the minimal Jacobi  operator $J_{0}$ has deficiency indices $(1,1)$.
 \end{corollary}
 
 Let us show that, even in the case $b_{n}=0$, assumption \e{eq:Gr6b} cannot be omitted in these assertions. Below we discuss a particular case of the examples considered  in \cite{Kost}, Example~2, and in \cite{J-M}, Lemma~2.3.
 
   \begin{example}\label{LCa}
Suppose that $b_{n}=0$ and that $a_{n}=n^p (1+c_{1}n^{-1})$ if $n$ is odd and $a_{n}=n^p (1+c_2n^{-1})$ if $n$ is even
for $p>1$ and sufficiently large $n$. Then
 \begin{equation}
\frac{a_{n}}{\sqrt{a_{n-1} a_{n+1}}}=1+ (-1)^n  \,\frac{c_{2}-c_{1}}{n} +O (\frac{1}{n^2})
\label{eq:Ber}\end{equation} 
as $n\to\infty$. Thus condition  \e{eq:nc} is now satisfied but \e{eq:Gr6b} fails if $c_{1} \neq c_{2}$.

It is proved in \cite{Kost, J-M} that the corresponding Jacobi  operator $J_{0}$ is essentially adjoint if $|c_{2}-c_{1}| \geq p-1$ while under the assumptions of \e{eq:nc}, \e{eq:Gr6b} the operators $J_{0}$ have deficiency indices $(1,1)$. 
This example shows also   that, without additional assumptions,  the Carleman condition  is not necessary for the essential self-adjointness of the operator $J_{0}$ even in the case $b_n=0$. Note that 
in view of \e{eq:Ber} condition $a_{n-1} a_{n+1}\leq a_{n}^2$ is now violated.
 \end{example}

 % We first prove that, for all $z\in\clos\Pi\setminus \{-1,1\}$, a solution $\{u_{n} (z)\}_{n=0}^\infty$ of equation \e{eq:A17} exists and then define the Jost solution  $\{f_{n} (z)\}_{n=0}^\infty$  of equation \e{eq:Jy} by formula \e{eq:A10}. 
 
 Recall that the polynomials   $P_{n}(z)$  are solutions of the Jacobi equation  \e{eq:Jy} satisfying the conditions  $P_{-1}(z)=0$,  $P_{0}(z)=1$.
 Relations \e{eq:LC}, \e{eq:LC1} and Theorem~\ref{LC} remain of course true in this case.  Moreover, the corresponding asymptotic relations in \e{eq:LC2} satisfy the relations $k_{-} (z)= \ov{k_{+}(\bar{z})}$ because $P_{n}(z)= \ov{P_{n}(\bar{z})}$. Set $\kappa(z)= |k_{+}(z)|$.
  Then
\[
k_{+} (z) = \kappa (z)e^{i\eta(z)},\q \eta(z)\in{\Bbb R}/2\pi{\Bbb Z}\q\mbox{and}\q   k_{-} (z) = \kappa (\bar{z})e^{-i\eta(\bar{z})}.
\]
Observe that the equalities $\kappa (z) = \kappa (\bar{z}) =0$ are impossible, since otherwise relation \e{eq:LC} for $P_{n}(z)$ would imply that $P_{n}(z)=0$ for all $n$ but $P_{0}=1$.

     The next result is a particular case of Theorem~\ref{LC}.  
     
     % if relations \e{eq:rt} and \e{eq:WW} are taken into account.

  \begin{theorem}\label{GPS}
    Let the assumptions of Theorem~\ref{GSS} be satisfied.  Then, for all $z\in{\Bbb C}$,  the sequence of the orthogonal polynomials $P_{n}( z )$  has  asymptotics
  \begin{equation}
P_{n}(  z )= a_{n} ^{-1/2}   \Big(\kappa (z) e^{-i\phi_{n}+i\eta(z)}+ \kappa (\bar{z}) e^{ i\phi_{n}- i\eta(\bar{z})}+ O( \varepsilon_{n})\Big), \q n\to\infty.
\label{eq:A22P}\end{equation} 
Asymptotics \e{eq:A22P} is  uniform in $z$ from compact subsets    of the complex plane $  \Bbb C $.
 \end{theorem}
 
  If $z=\lambda\in {\Bbb R}$, then relation  \e{eq:LC} for $P_{n}(  z )$ reduces to
   \begin{equation}
P_{n}(\lambda)= 2 \Re \big(k_{+}(\lambda) f_{n}(\lambda) \big),
\label{eq:WW2}\end{equation}
and \e{eq:A22P}  yields the following result.

   \begin{corollary}\label{GPSc}
 If $\lambda\in {\Bbb R}$, then
  \begin{equation}
P_{n}(  \lambda )= 2 a_{n} ^{-1/2}   \Big(\kappa (\lambda) \cos \big( \phi_{n}-\eta(\lambda)\big) + O( \varepsilon_{n})\Big), \q n\to\infty.
\label{eq:A22r}\end{equation} 
Asymptotics \e{eq:A22r} is  uniform in $\lambda$ from compact subsets    of the line $  \Bbb R $.
 \end{corollary}
 
 Observe that $\kappa (\lambda) \neq 0$, since otherwise \e{eq:WW2} would imply that $P_{n}(\lambda)=0$ for all $n$ but $P_{0}=1$.  
 
 Next we discuss specially the case   $\Im z \neq 0$. Multiplying equation \e{eq:RR} for $P_n (z)$ by  its complex conjugate $\bar{P}_n (z)$ and taking the sum over $n=0,1,\ldots, N$, we find that
  \begin{equation}
 \sum_{n=0}^{N}  a_{n-1} P_{n-1}  (z) \bar{P}_n (z) + \sum_{n=0}^{N} a_{n} P_{n+1} (z) \bar{P}_n (z) + \sum_{n=0}^{N} b_{n} |P_{n} (z)|^{2}  = z  \sum_{n=0}^{N}  |P_{n} (z)|^{2} . 
\label{eq:Jxy}\end{equation}
Since $P_{-1}  (z)=0$,
the first sum on the left equals
\[
 \sum_{n=0}^{N}  a_{n} P_{n}  (z) \bar{P}_{n+1} (z) -  a_N P_{N}  (z)\bar{P}_{N+1} (z).
 \]
 Therefore taking the imaginary part of \e{eq:Jxy}, we see that
   \begin{equation}
a_N \Im   \big(P_{N+1}  (z)\bar{P}_N (z) \Big) = \Im z \sum_{n=0}^{N}  |P_{n} (z)|^{2} . 
\label{eq:Jxy1}\end{equation}
According to \e{eq:A22P}    the left-hand side here equals
 \begin{equation}
 \varkappa_N^{-1} \Im  \Big(  \big(\kappa (z) e^{-i\phi_{N+1}+i\eta(z)}+ \kappa (\bar{z}) e^{ i\phi_{N+1}- i\eta(\bar{z})} \big)
 \big(\kappa (z) e^{i\phi_{N }-i\eta(z)}+  \kappa (\bar{z}) e^{- i\phi_{N}+ i\eta(\bar{z})} \big) \Big)+ o(1)
\label{eq:Jxy2}\end{equation} 
with $ \varkappa_N$ defined by  \e{eq:Gr4K}.
Observe that the sum 
  \begin{multline*}
\kappa (z) e^{-i\phi_{N+1}+i\eta(z)} \kappa(\bar{z}) e^{- i\phi_{N}+ i\eta(\bar{z})}
+\kappa(\bar{z}) e^{ i\phi_{N+1}- i\eta(\bar{z})} \kappa(z) e^{i\phi_{N }-i\eta(z)}\\=2 \kappa (z) \kappa(\bar{z})  \cos\big( \phi_{N }+\phi_{N+1} -\eta(z)-\eta(\bar{z}) \big)
 \end{multline*}
is real. Therefore taking relations   and \e{eq:Gr6a},  \e{eq:Gs3}  and  \e{eq:GRR} into account, we see that  expression \e{eq:Jxy2} equals
  \[
\varkappa_N^{-1}   \Im    \big(\kappa (z)^{2} e^{-i\theta_{N }}+ \kappa (\bar{z})^{2} e^{i\theta_{N } }\big)+ o(1)
=
\varkappa_{\infty}^{-1}    \big(  \kappa (\bar{z})^{2} -\kappa (z)^{2}   \big) \sin \theta_\infty + o(1)
\]
where $\sin \theta_\infty=\sqrt{1-\beta_{\infty}^2}$.
Replacing the left-hand side of
\e{eq:Jxy1} by this  expression, we obtain the following result.

   \begin{theorem}\label{zz}
 Under the assumptions of Theorem~\ref{GSS} the  identity
 \[
\kappa (\bar{z})^{2} - \kappa (z)^{2} =\Im z \: \varkappa_{\infty}  (1-\beta_{\infty}^2)^{-1/2}  \sum_{n=0}^{\infty}  |P_{n} (z)|^{2}
  \]
  holds.
 \end{theorem}
 
   \begin{corollary}\label{zz1}
   If $\Im z>0$, then
     $ \kappa (z)< \kappa (\bar{z})$. Equivalently,
  $ \kappa (z)> \kappa (\bar{z})$   if $\Im z<0$.
      \end{corollary}

 % \subsection{Dominating diagonal elements}
 
       \subsection{Large diagonal elements}
       
 Here we consider the case      $|\beta_{\infty}|>1$. Choose an arbitrary $z\in{\Bbb C}$. By Theorem~\ref{GSS+}, the sequence $f_{n}( z )$ defined by equality \e{eq:JostG} satisfies equation \e{eq:Jy},  and it has exponentially decaying asymptotics \e{eq:A22G+} where the phase $\varphi_{n}$ is given by \e{eq:Grf1}. In particular, $f_{n}(z)\neq 0$ for sufficiently large $n$, say, $n\geq n_{0}= n_{0}  (z)$.
 Now we have to construct a solution $g_{n}( z )$ of \e{eq:Jy} linearly independent with $f_{n}( z )$.
 We define it by the formula \e{eq:GEg}, that is,
    \begin{equation} 
g_{n}(z) = f_{n}(z) G_{n}(z)
\label{eq:GH}\end{equation}
 where
   \begin{equation} 
 G_{n} (z)=\sum_{m=n_{0}}^n (a_{m-1}f_{m-1}(z) f_{m}(z))^{-1},\q n \geq n_{0} .
\label{eq:GE+}\end{equation}

%where $n_{0}$ is chosen in such a way that ${\sf f}_{m-1}(z)\neq 0$ for all 

\begin{theorem}\label{GE+}
Suppose that a  sequence $ f (z)=\{f_{n}(z)  \}$  satisfies the Jacobi equation \e{eq:Jy}.
Then   the sequence $ g (z)=\{g_{n}(z)  \}$ defined by formulas  \e{eq:GH} and \e{eq:GE+}
 satisfies the same equation  and
   the Wronskian $\{ f(z), g(z)\}=1$. In particular, 
   the solutions $f(z)$ and $g (z)$ are linearly independent.
  \end{theorem}

\begin{pf}
First, we check equation \e{eq:Jy} for $g_{n}$.  Observe that
 \begin{multline*}
a_{n-1} f_{n-1} G_{n-1}+ (b_{n} -z) f_{n} G_{n} +a_{n} f_{n+1}G_{n+1}
\\
=\big(a_{n-1} f_{n-1}  + (b_{n} -z) f_{n}  +a_{n} f_{n+1}\big) G_{n}
+ a_{n-1} f_{n-1} ( G_{n-1}- G_{n}) + a_{n} f_{n+1}(G_{n+1}-G_{n}) .
\end{multline*}
The first term here is zero because    equation \e{eq:Jy} is true for the  sequence $f_{n}$.  According to definition \e{eq:GE+}
   \begin{equation} 
G_{n+1}=G_n+ (a_{n} f_{n} f_{n+1})^{-1},
\label{eq:GEx}\end{equation}
so that the second and third terms equal $- f_{n}^{-1}$ and $f_{n}^{-1}$, respectively. This proves equation  \e{eq:Jy} for $g_{n}$.

It also follows from \e{eq:GEx}
that the Wronskian \e{eq:Wr}   equals
  \[
\{ f(z), g(z)\}= a_{n}f_{n}(z)f_{n+1}(z)(G_{n+1}(z)- G_{n}(z))=1,
\]
 whence the solutions $f(z)$ and $g(z)$ are linearly independent.
\end{pf}

 % of $f (z)=\{f_{n}(z)\}$ and $g (z)=\{f_{n}(z) G_{n}(z)\}$
 
The following statement shows that   the solution $g_{n}  (z)$ of equation \e{eq:Jy}  exponentially grows as $n\to\infty$. 
For a proof, we will have to integrate by parts in \e{eq:GE+}. It is convenient to start with a simple technical assertion.
Recall that $\sigma_{m} =\z_{m-1} \z_{m}=e^{-\vartheta_{m-1}-\vartheta_{m}} $, the numbers $\varkappa_{n}$ were defined in  formula \e{eq:Gr4K} and
 the sequence $u_n (z)$ was constructed in Theorem~\ref{GS3}.
%$\z_{n}$ is defined by formula \e{eq:Gs3+} and $\varphi_{n}$ is defined in \e{eq:A22G+}.    % We also set
%\begin{equation}
%\tau_{m}  = \z_{m-1}\z_{m} (1-\z_{m-1}\z_{m})^{-1}.
%\label{eq:tt}\end{equation}

 \begin{lemma}\label{GEr}
Let
\begin{equation}
 v_{m} (z)=  \sgn\beta_{\infty} \varkappa_{m-1} (\sigma_{m}^{-1}  -1)^{-1} (u_{m-1} (z) u_{m} (z))^{-1}  
 \label{eq:hv}\end{equation}
and
\begin{equation} 
t_{m}=e^{\varphi_{m-1}+\varphi_{m}} .
\label{eq:hvv}\end{equation}
Then 
 \begin{equation}
(a_{m-1} f_{m-1} (z)f_{m}(z))^{-1}=    v_{m} (z) t'_{m}.
\label{eq:hv3}\end{equation}
 \end{lemma}
 
 \begin{pf}
 It follows from formulas \e{eq:Gr1+} and \e{eq:JostG} that
 \[
 a_{m-1} f_{m-1} f_{m} =\sgn\beta_{\infty} \frac{a_{m-1}}{\sqrt{a_{m-1} a_{m}}}   u_{m-1} u_{m} t_{m}^{-1},
 \]
 whence
 \[
(a_{m-1} f_{m-1} f_{m})^{-1}= \sgn\beta_{\infty} \varkappa_{m-1} (u_{m-1} u_{m})^{-1} t_{m} .
\]
By definitions \e{eq:GRR} and \e{eq:hvv}, we have
\[
t_{m}'= t_{m}  (\sigma_{m}^{-1}-1).
\]
Putting together the last two formulas with the definition \e{eq:hv}  of $v_{m} $, we arrive at \e{eq:hv3}.
 \end{pf}
 
% \begin{equation}
%  t_{m}= \tau_{m} t_{m}' 
%\label{eq:hv1}\end{equation}
%and, by virtue of \e{eq:Gr4K}, 
%\[(a_{m-1} {\sf f}_{m-1} {\sf f}_{m})^{-1}=  (\varkappa_{m-1} {\sf u}_{m-1} {\sf u}_{m})^{-1}t_{m}=  v_{m} t'_{m}.\]

Using \e{eq:hv3} and 
integrating  by parts   in \e{eq:GE+}, we  find that
  \begin{align}
G_{n} =\sum_{m=n_{0}}^n v_{m}  t_{m}' = v_{n}  t_{n+1} -  v_{n_{0}-1}  t_{n_{0}} - \wt{G_{n}} \q \mbox{where} \q \wt{G_{n}}
=   \sum_{m= n_{0}}^n v_{m-1}'   t_{m}.
\label{eq:GE1}\end{align}
We will see that the asymptotics of sequence 
\e{eq:GH} as $n\to \infty$ is determined by the first  term in the right-hand side of \e{eq:GE1}. Let us calculate it.

 \begin{lemma}\label{GErY}
The asymptotic relation
     \begin{equation}
\lim_{n\to\infty} \sqrt{a_{n}} e^{-\varphi_{n}} (\sgn \beta_{\infty})^n  f_{n} (z)  v_{n} (z) t_{n+1} =  \frac{ \sgn \beta_{\infty} }{2\sqrt{ \beta_{\infty}^2 -1} }   \varkappa_{\infty}
\label{eq:gar1}\end{equation}
holds.
 \end{lemma}

  \begin{pf}
  It follows from Theorem~\ref{GSS+} that
  \[
 \sqrt{a_{n}}  (\sgn \beta_{\infty})^n  f_{n} (z)=  e^{-\varphi_{n}} (1+ o(1)).
  \]
Using   definition \e{eq:hv}  of $v_n $, Lemma~\ref{GK} and Theorem~\ref{GS3} we find that
\[
 v_n (z)=    \sgn \beta_{\infty} \varkappa_\infty( \z_\infty ^{-2} -1)^{-1} (1+ o(1))  .
\]
%label{eq:hvx}\end{equation}
Since according to   \e{eq:hvv}  
    \[
t_{n+1}=e^{2\varphi_{n}+\vartheta_{n}} =e^{2\varphi_{n} }  \z_\infty^{-1} (1+ o(1)),
\]
%label{eq:hv0}\end{equation}
the limit \e{eq:gar1} equals $ \sgn \beta_{\infty}  \varkappa_{\infty}  \z_\infty^{-1}( \z_\infty^{-2} -1)^{-1}$ which in view of \e{eq:zz+} coincides with the right-hand side of \e{eq:gar1}.
     \end{pf}
     
     Next we show that the remainder $\wt{G}_{n} (z)$ in \e{eq:GE1} is negligible. We use the following technical assertions.
     
        \begin{lemma}\label{GE1X}
The sequence $v_{n} (z)$   defined by \e{eq:hv} satisfies the condition  $\{v_{n}' \}\in \ell^1 ({\Bbb Z}_{+})$.
 \end{lemma}
 
     \begin{pf}
     According to \e{eq:Gs3+} we have
     \[
|\z_{n}|= (  |  \beta_n |+ \sqrt{  \beta_n^2-1})^{-1} <|\beta_{\infty}|^{-1}<1
\]
%label{eq:G33+}\end{equation}
for sufficiently large $n$. Moreover, it follows from \e{eq:Gs2} that $\{\z_{n}' \}\in \ell^1 ({\Bbb Z}_{+})$ if $\{\beta_{n}' \}\in \ell^1 ({\Bbb Z}_{+})$ whence
\[
\big(  ( \sigma_{n}^{-1} -1)^{-1}\big)' \in \ell^1 ({\Bbb Z}_{+}).
\]
 Putting together Theorem~\ref{GS3}  and Lemma~\ref{GS'} we see that
\[
(u_{n-1}  ^{-1} u_{n}  ^{-1})' \in \ell^1 ({\Bbb Z}_{+}).
\]
Now the inclusion $\{v_{n}' \}\in \ell^1 ({\Bbb Z}_{+})$ is a direct consequence of definition \e{eq:hv}.
         \end{pf}
     
      \begin{lemma}\label{GErX}
Let $\wt{G}_{n} (z)$ be defined in formula \e{eq:GE1}. Then
     \begin{equation}
\lim_{n\to\infty} \sqrt{a_{n}} e^{-\varphi_{n}} f_{n}  (z)  \wt{G}_{n} (z)=  0.
\label{eq:gar}\end{equation}
 \end{lemma}
 
    \begin{pf}
Using Theorem~\ref{GSS+} and definition \e{eq:hvv}, we see that
\begin{multline*}
\sqrt{a_{n}} e^{-\varphi_{n}} |f_{n}\wt{G}_{n} |\leq C e^{-2\varphi_{n}} \sum_{m=n_{0}}^n  |v_{m-1}'|e^{\varphi_{m-1}+\varphi_{m}}
\\
 \leq C_{1} e^{-\varphi_{n}}\sum_{n_{0}\leq m <[n/2]} |v_{m-1}'| + C_{1}  \sum_{[n/2]\leq m \leq n} |v_{m-1}'|.
%\label{eq:Gt3}
\end{multline*}
Since $\varphi_{n}\to\infty$ as $n\to\infty$, relation  \e{eq:gar} follows from
 Lemma~\ref{GE1X}.
    \end{pf}
     
   Let us come back to the representations \e{eq:GH} and    \e{eq:GE1}.
       Putting together Lemmas~\ref{GErY}, \ref{GErX} and taking into account that the term $ v_{n_{0}-1}  t_{n_{0}}$   is negligible, we obtain the asymptotics of $g_{n} (z)$.

 \begin{theorem}\label{GEe}
% Let $|\beta_{\infty}| >1$, condition   \e{eq:asqr1} be satisfied.
   Under the assumptions of Theorem~\ref{GSS+} the solution \e{eq:GEg} of the Jacobi equation  \e{eq:Jy} satisfies for all $z\in{\Bbb C}$ the relation
     \[
\lim_{n\to\infty} \sqrt{a_{n}} e^{-\varphi_{n}}  (\sgn \beta_{\infty})^{n } g_{n} (z)=   \frac{\sgn \beta_{\infty}  }{2\sqrt{ \beta_{\infty}^2 -1} }\varkappa_{\infty}  .
\]
%label{eq:gas+}\end{equation}
 \end{theorem}
 
 In view of \e{eq:xz+}, \e{eq:zz1}, $\{ g_{n} (z)\}\in \ell^2 ({\Bbb Z}_{+})$ if and only if
     \begin{equation}
  \sum_{n =0}^{\infty} a_{n}^{-1} \big(|\beta_{\infty}|+\sqrt{\beta_{\infty}^2 -1} \big)^{2n}<\infty . 
 \label{eq:nc-}\end{equation} 
 In this case all solutions of equation \e{eq:Jy} are in $\ell^2 ({\Bbb Z}_{+})$. Therefore we can state
 
  \begin{corollary}\label{GEc}
   Under the assumptions of Theorem~\ref{GSS+} the minimal Jacobi operator $J_{0}$ is essentially self-adjoint if and only if      condition \e{eq:nc+} is satisfied. Otherwise, that  is under assumption \e{eq:nc-},  the operator $J_{0}$  has deficiency indices $(1,1)$.
 \end{corollary}
 
 In particular, the  operator $J_{0}$ is essentially self-adjoint  if $a_{n}\leq cn^p$  for some $c> 0$, $p<\infty$ and  all $n\geq 1$. 
 
  \begin{example}\label{GEd}
  Suppose that $a_{n}= c_{0} x^n (1+o(1))$ where $c_{0} >0$ and $x>1$ as $n\to\infty$. Then condition \e{eq:nc+} is satisfied if and only if
  \[
\sqrt{x}  \leq |\beta_{\infty}|+\sqrt{\beta_{\infty}^2 -1}  .
  \]
 \end{example}
 
  \begin{example}\label{GEex}
  Suppose that $a_{n}= c_{0} x^{n^{p}} (1+o(1))$ where $c_{0} >0$   as $n\to\infty$. Then condition \e{eq:nc-} is satisfied  for all $p>1$,  $x>1$ and all $|\beta_{\infty}|>1$.
 \end{example}
 
    Set $ P(z)=\{ P_{n}(z)\}_{n=-1}^\infty$, $ f(z)=\{ f_{n}(z)\}_{n=-1}^\infty$, and define the Jost function by the equality
    \begin{equation}
\Omega (z):=  \{ P(z), f(z)\} = - 2^{-1}f_{-1}(z),
\label{eq:WRG}\end{equation}
where the first formula \e{eq:Wr1}   has been used. 
 By Theorem~\ref{GE+}, the Wronskian $\{ f(z),g(z)\}=1$ so that
 \[
P_{n} (z)= \omega(z)f_{n} (z)-\Omega (z)g_{n} (z)
\]
with $\omega(z)= \{ P(z),g(z)\}$.  Obviously,  
$\omega (z)\neq 0$ if $\Omega(z)= 0$. Therefore Theorems~\ref{GSS+}   and \ref{GEe} imply the following result.

\begin{theorem}\label{GE1}
   Under the assumptions of Theorem~\ref{GSS+} the relation 
    \begin{equation}
  P_{n}(z)=-    \varkappa_{\infty} \Omega(z)  \frac{ (\sgn \beta_{\infty})^{n +1} }{2\sqrt{\beta^2_{\infty}-1}}  \frac{e^{\varphi_{n}}  }{  \sqrt{a_{n}}} (1+ o(1)), \q n\to\infty,
\label{eq:GEGE}\end{equation}
is true for all  $z\in{\Bbb C} $  
with convergence uniform on compact subsets of $z\in{\Bbb C} $.
Moreover, if $\Omega(z)=0$, then
 \[
\lim_{n\to\infty} \sqrt{a_{n}} e^{\varphi_{n}} (\sgn \beta_{\infty})^{n } P_{n}(z)=\{ P(z),g(z)\} \neq 0.
\]
%label{eq:GEGEx}\end{equation}
 \end{theorem}
 
 If condition \e{eq:nc-} is satisfied, 
     then the 
closure $J $ of  the minimal Jacobi operator $J_{0}$ is   self-adjoint. In this case $\Omega(z)=0$ if and only if $z$ is an eigenvalue of $J$.
The resolvent  of the operator $J$ can be constructed by  the standard (cf. Lemma~2.6 in \cite{Y/LD}
or Lemma~5.1 in \cite{JLR}) formulas. Recall that $e_{n}$, $n\in {\Bbb Z}_{+}$, is the canonical basis in $\ell^2 ({\Bbb Z}_{+})$, $f_{n} (z)$ is the Jost solution of the   equation \e{eq:Jy} constructed in Theorem~\ref{GSS+} and 
   $\Omega(z)$ is the Wronskian \e{eq:WRG}.

   \begin{proposition}\label{resolvent}
Under the    assumptions of Theorem~\ref{GSS+} suppose also that condition \e{eq:nc+} holds. Then the resolvent $(J -z )^{-1}$ of the Jacobi operator  $J$  satisfies the relations
   \begin{equation}
((J -z )^{-1} e_{n}, e_{m})= \Omega(z)^{-1} P_{n} (z) f_{m}(z),\q \Im z \neq 0, 
\label{eq:RRe}\end{equation}
if $n\leq m$ and $((J -z )^{-1} e_{n}, e_{m})=((J -z )^{-1} e_m, e_n)$.
 \end{proposition}

Since in view of Theorem~\ref{GSS+}, $f_{n} (z)$ and hence   $\Omega(z)$ are entire functions of $z\in{\Bbb C}$,  we can state

\begin{corollary}\label{rdiscr}
The spectrum of the  operator  $J$  is discrete, and its  eigenvalues $\lambda_{1}, \cdots, \lambda_{k}, \ldots$   are given by the equation  $\Omega(\lambda_{k}) =0$.  The resolvent $(J -z )^{-1} $ is an analytic function  of  $z\in{\Bbb C}$ with poles in the points $\lambda_{1}, \cdots, \lambda_{k}, \ldots$.
 \end{corollary}
 
  A much more general result of this type is stated below as Proposition~\ref{discr}.
 
 In view of formula \e{eq:WRG} and equation \e{eq:Jy} where $n=0$ for the Jost solution $f_{n} (z)$, the equation for eigenvalues of $J$ can be also written  as
  \[
       (b_{0}-\lambda_{k})f_0(\lambda_{k})+ a_{0} f_1(\lambda_{k})=0.
       \]
       It also  follows from representation \e{eq:RRe} where $n=m=0$ and the first formula for \e{eq:Wr1} that the spectral measure of $J$ is given by the standard formula
         \[
\rho(\{\lambda_{k}\})= 2\frac{f_{0}(\lambda_{k})} {\dot{f}_{-1}(\lambda_{k})}  
\]
where $\dot{f}_{-1}(\lambda )$ is the derivative of  $f_{-1}(\lambda )$ in $\lambda$. Alternatively, we have
       \begin{equation}
\rho(\{\lambda_{k}\})= \big(\sum_{n=0}^\infty P_{n}(\lambda_{k})^2\big)^{-1}. 
\label{eq:meas}\end{equation}
A proof of this relation can be found, for example, in \cite{Simon}, Theorem~4.11.
       
     \section{ Discussion} 
     
        \subsection{Operators with discrete spectrum} 
        
        We start with a very general result which, in particular, applies to Jacobi operators. An assertion below is essentially known, and we give it mainly for completeness of our presentation. The operator $J$ will be defined via its quadratic form
          \begin{equation}
(J u,u)=  \sum_{n=0}^\infty b_{n}|u_{n}|^2 + 2 \Re \sum_{n=0}^\infty a_{n}u_{n} \bar{u}_{n+1}
\label{eq:QF}\end{equation}
where $a_{n}$ are complex and $b_{n}$ are real numbers.

           \begin{proposition}\label{discr}
                    Suppose that $b_{n}\to \infty$ $($or $b_{n}\to -\infty)$ and 
that
   \begin{equation}
|a_{n-1}|+ |a_{n }|\leq \epsilon |b_{n}|
\label{eq:QF1}\end{equation}
for some $\epsilon<1$ and all sufficiently large $n$.
Then the spectrum of the operator $J$ is discrete and is semi-bounded from below $($resp.,  from above$)$.
 \end{proposition}
           
         \begin{pf}      
          Since finite-rank perturbations cannot change the discreteness of  spectrum, we can suppose that estimate 
   \e{eq:QF1} is true for all $n$.
     Let $B$ be the operator corresponding to the form 
           \begin{equation}
(B u,u)=  \sum_{n=0}^\infty b_{n}|u_{n}|^2 .
\label{eq:QF2}\end{equation}
Evidently, its spectrum is discrete and is semi-bounded from below (resp.,  from above).  For a proof of the same result for the operator $J$, it suffices to check that the forms   \e{eq:QF} and   \e{eq:QF2} are equivalent or that
     \begin{equation}
  2 \big|\Re \sum_{n=0}^\infty a_{n}u_{n} \bar{u}_{n+1} \big| \leq \epsilon \sum_{n=0}^\infty| b_{n}| |u_{n}|^2 
\label{eq:QF3}\end{equation}
for some $\epsilon<1$. Let us use the following obvious inequality
     \[
  2 \big|  \sum_{n=0}^\infty a_{n}u_{n} \bar{u}_{n+1} \big| \leq 2\sqrt{\sum_{n=0}^\infty | a_{n}| |u_{n}|^2 \sum_{n=0}^\infty | a_{n}|
  |u_{n+1}|^2} \leq  \sum_{n=0}^\infty (| a_{n-1}| + | a_{n}|)
  |u_{n}|^2 .
\]
%label{eq:QF4}\end{equation}
Therefore estimate \e{eq:QF3} is a direct consequence of the condition \e{eq:QF1}.
       \end{pf}       
       
  Proposition~\ref{discr} holds, in particular,  for semi-bounded Jacobi operators when $a_{n} >0$. It applies directly to the Friedrichs' extension of the operator $J_{0}$, but its conclusion remains true for all extensions $J$ of  $J_{0}$ because the deficiency indices of $J_{0}$ are finite.
   In particular, condition \e{eq:QF1} is satisfied with any $\epsilon \in(|\beta_{\infty}|^{-1},1)$ under the assumptions of Theorem~\ref{GSS+}. Thus,    we recover the result of Corollary~\ref{rdiscr}. 
   Clearly,  the corresponding operator $J $  is semi-bounded from below (from above) if $\beta_{\infty}<-1$ (if $\beta_{\infty} >1$, respectively).

       \subsection{The Carleman case } 
       
        In this subsection we  still consider the case $a_{n}  \to\infty$ but assume that  the Carleman condition  \e{eq:Carl} holds. Other assumptions on $a_{n}$ and $b_{n}$ are essentially the same as in the main part of this paper. In particular, we suppose that  condition  \e{eq:Gr} is satisfied.
         Recall that under the Carleman condition   the minimal Jacobi operator $J_{0}$ is essentially self-adjoint.    The results stated here  will be published elsewhere, but some of them are close to  the papers \cite{Janas, Apt, Sw-Tr}.
       
% In this subsection we  still consider the case $a_{n}  \to\infty$ but assume that
%  \begin{equation}\sum_{n=0}^\infty \alpha_{n} = \infty\label{eq:asqr1}\end{equation}
%where the numbers $\alpha_{n} $ are defined in \e{eq:aabb}. This  assumption is slightly stronger than
 %   the Carleman condition  \e{eq:Carl}, and \e{eq:asqr1} is equivalent to  \e{eq:Carl} if the limit \e{eq:Gr6a} exists.  It is convenient to assume that $\varkappa_{\infty}  =1$.  Recall that under the Carleman condition   the minimal Jacobi operator $J_{0}$ is essentially self-adjoint. Assumption   \e{eq:Gr} on the coefficients $b_{n}$ is the same as in the non-Carleman case.     The results stated here  will be published elsewhere, but some of them are close to  the papers \cite{Janas, Apt, Sw-Tr}.
 
 % , but we restrict ourselves to the case $|\beta_{\infty}|<1$
       
 Now the term $2z\alpha_{n}$ in the right-hand side of \e{eq:Grr1}  cannot be neglected and the numbers $\z_{n}$ have to be defined as approximate solutions of the equation
  \[
 \z_{n }+  \z_{n}^{-1}=  2 \beta_{n} + 2z \alpha_{n}  .
\]
%label{eq:ANS1}\end{equation}  
  For simplicity of our discussion,  we assume here that
        \begin{equation}
\sum_{n=0}^\infty   a_{n}^{-3}  (1+ |b_{n}|)<\infty.
\label{eq:D}\end{equation} 
Conditions \e{eq:Carl} and \e{eq:D} admit a growth of the off-diagonal coefficients $a_{n}$ as $n^p$ for $p\in (1/2, 1]$ and even for $p\in (1/3, 1]$ if $b_{n}=0$.

 Suppose first that $|\beta_{\infty}|<1$.
  As before, the Ansatz $Q_{n}  (z)$ is defined by formula \e{eq:ANS}, but, instead of  \e{eq:Gs3},
   we now set
        \begin{equation}
\z_{n}=   ( \beta_n- i\sqrt{ 1- \beta_n^2}) \exp\big( iz\frac{\alpha_{n}}{\sqrt{1-\beta_{n}^2}}\big).
\label{eq:GS3}\end{equation}
An easy calculation shows that the corresponding remainder \e{eq:Grr1} is in $\ell^1 ({\Bbb Z}_{+})$. Therefore repeating the arguments of  Sect.~2 and 3, we find that equation \e{eq:Jy} has a solution 
(the Jost solution) $f_{n} (z)$   distinguished by the asymptotics 
   \begin{equation}
f_{n} (z)=a_{n}^{-1/2}e^{- i \phi_{n} + iz \psi_{n} }(1+ o(1)), \q \Im z \geq 0, \q n\to \infty,
\label{eq:D1}\end{equation} 
where the phase $\phi_{n}$ is defined by formula \e{eq:Grf} and
 \begin{equation}
\psi_{n}= \sum_{m=0}^{n-1} \ \frac{ {\alpha}_m }{\sqrt{|1-\beta_m^{2}|}}.
\label{eq:psi}\end{equation} 
It follows from condition \e{eq:D} that $\psi_{n}=O (n^{2/3})$ as $n\to\infty$ so  that by virtue of \e{eq:zz1} the phase $\psi_{n}$ is negligible compared to $\phi_{n}$.
We set $f_{n}(z)=\ov{f_{n}(\bar{z})}$ for $\Im z \leq 0$.
The functions $f_{n} (z)$ are analytic in the complex plane with the cut along $\Bbb R$ and are continuous up the cut. It is easy to see that $f_{n} (z)\in \ell^2 ({\Bbb Z}_{+})$ for $\Im z \neq 0$.
% Indeed, it  follows from \e{eq:D1} that
%   \[ | f_{n} (z) | =a_{n}^{-1/2}\exp \big(- | \Im z |\,  \psi_{n} \big) (1+ o(1))   \]
%  and therefore it suffices to use an elementary inequality
 %      \begin{equation}
 %  \sum_{n=1}^\infty \alpha_{n}\exp\Big(-  \sum_{m=0}^{n-1} \alpha_m\Big)<\infty.
 %  \label{eq:asqrf2}\end{equation} 
%   For a proof of the convergence of this series, we can replace $\alpha_{n}$ by $1-e^{-\alpha_{n}}$ and then observe that
%      \begin{multline*}
%      \sum_{n=1}^N (1-e^{-\alpha_{n}})\exp\Big(-  \sum_{m=0}^{n-1} \alpha_m\Big)\\
%      =  \sum_{n=1}^N  \exp\Big(-  \sum_{m=0}^{n-1} \alpha_m\Big)- \sum_{n=2}^{N+1}  \exp\Big(-  \sum_{m=0}^{n-1} \alpha_m\Big)
  %    =e^{-\alpha_{0}}- \exp\Big(-  \sum_{m=0}^N \alpha_m\Big)< e^{-\alpha_{0}}.
  %  \end{multline*}
%    By virtue of \e{eq:asqrf} estimate \e{eq:asqrf2} implies \e{eq:asqrf1}. 
   
   If $z =\lambda\in {\Bbb R}$, we have two Jost solutions $f_{n}(\lambda+i0)$ and $f_{n}(\lambda-i0)=\ov{f_{n}(\lambda+i0)}$ of equation \e{eq:Jy}. Their Wronskian is the same as \e{eq:A2CG1} (where $\varkappa_{\infty}=1$):
     \[
\{ f (\lambda+i0),  f(\lambda-i0)\} 
= 2 i   \sqrt{1- \beta^2_{\infty}} \neq 0.
 \]
 Therefore the   polynomials $P_{n}  (\lambda)$ are linear combinations of $f(\lambda+i0)$ and $  f (\lambda-i0)$:
  \begin{equation}
P_{n} (\lambda)=\frac{  \Omega (\lambda- i0)  f_{n} (\lambda+i0) - \Omega(\lambda + i0)  f_{n} (\lambda- i0)  }{ 2 i   \sqrt{1- \beta^2_{\infty}}},   \q n=0,1,2, \ldots, \q \lambda\in \Bbb R,
\label{eq:HH4L}\end{equation}
  where $\Omega (z)$ is defined by formula \e{eq:WRG}. 
 Note that $\Omega (\lambda\pm i0)\neq 0$ for all $\lambda\in \Bbb R$ according to \e{eq:HH4L}.
 It now follows from \e{eq:D1} for $z =\lambda\in {\Bbb R}$ that   the  polynomials $P_{n} (\lambda)$ have asymptotics 
   \begin{equation}
 P_{n} (\lambda)= - a_{n}^{-1/2} \Big(   | \Omega ( \lambda+i0)| (1- \beta^2_{\infty})^{-1/2}     \sin (\phi_{n}  -\lambda\psi_{n} +\arg \Omega(\lambda+i 0) ) + o(1) \Big) , \q n\to\infty.
\label{eq:Sz}\end{equation}
%where
% \begin{equation}
% \Phi_{n}(\lambda )= -\phi_n  +   \lambda \psi_{n}   .
%\label{eq:Jo1}\end{equation}

  If $\Im z \neq 0$, then by virtue of Theorem~\ref{GE+} a solution  $g_{n}(z)$ of the Jacobi equation \e{eq:Jy} linear independent with $f_{n}(z)$
    can be constructed by   explicit formula \e{eq:GEg}. Similarly to Theorem~\ref{GEe}, it can be checked that
   %     It exponentially grows as $n\to\infty$, namely 
  \[
   g_{n}(z)=\frac{1}{2i \sqrt{1 - \beta^2_{\infty} } }  \frac{1}{ \sqrt{a_{n}}  }e^{i \phi_{n} - iz \psi_{n}}   (1+o(1)), \q \Im z>0, 
   \q n\to\infty,
\]
%label{eq:gas}\end{equation}
so that $   g_{n}(z)$ grows as $n\to\infty$ faster than any power of $n$.
Since the Wronskian $\{f(z), g(z)\} =1$, the  asymptotics   of the orthonormal polynomials can be easily derived from this result:
  \begin{equation}
    P_{n}(z)=\frac{i \Omega(z)}{2 \sqrt{1 - \beta^2_{\infty} } } 
 \frac{1}{ \sqrt{a_{n}}  }  e^{i \phi_{n} - iz \psi_{n}}  (1+o(1)), \q \Im z>0,    \q n\to\infty.
\label{eq:Gas}\end{equation}
   
     Under assumption  \e{eq:Carl}   the Jacobi operator $J_{0}$ is essentially self-adjoint (see, e.g., the book  \cite{AKH}).  
 The resolvent $ (J-z)^{-1}$ of $J=\clos J_{0}$ is given by the general formula \e{eq:RRe}. 
 Since the Jost solutions $f_{n} (z)$ are continuous functions of $z$ up to the real axis and 
 $\Omega (\lambda\pm i0)\neq 0$ for   $\lambda\in \Bbb R$, the spectrum of  the Jacobi operator $J$ is absolutely continuous. Moreover,
 using   relation  \e{eq:HH4L}, it is easy to deduce from \e{eq:RRe} a representation for the spectral family $d E (\lambda)$ of the operator $J$:
   \[
 \frac{d(E (\lambda)e_n, e_m)} {d\lambda}= \pi^{-1}\sqrt{1- \beta^2_{\infty}} \, | \Omega (\lambda+i0) |^{-2} P_{n} (\lambda) P_m (\lambda)   .
\]
  In particular, for  the spectral measure, we have the expression
  \[
  d\rho(\lambda)=  \pi^{-1}\sqrt{1- \beta^2_{\infty}} \, | \Omega (\lambda+i0) |^{-2} d\lambda.
  \]
  It follows that the spectrum of   $J$   coincides with the whole real axis.

  Formulas \e{eq:Sz} and \e{eq:Gas} are consistent with the classical asymptotic expressions for the Hermite polynomials  when $a_{n }=\sqrt{(n+1)/2}$ and $b_{n}= 0$ (see, e.g., Theorems~8.22.6 and 8.22.7 in the G.~Szeg\H{o}'s book \cite{Sz}). As far as previous general results for $\Im z\neq 0$ are concerned, we are aware only of the paper \cite{Rah} where
    an asymptotics       of $| P_{n}(z)|$ as $n\to\infty$   was found   in terms of   a behavior of the corresponding absolutely continuous spectral measure  for $|\lambda|\to\infty$.

      The cases $|\beta_{\infty}|> 1$   and $|\beta_{\infty}|< 1$ are technically rather similar, but the asymptotic behavior of orthogonal polynomials and spectral properties of the Jacobi operators are quite different in these cases.  If $|\beta_{\infty}|> 1$,  we set (cf. \e{eq:GS3})
    \[
\z_{n}=  \sgn\beta_{n} ( |\beta_n | -  \sqrt{ \beta_n^2-1})  \exp\big( -z\frac{\alpha_{n}}{\sqrt{\beta_{n}^2-1}}\big).
\]
%label{eq:GS3+}\end{equation}
Then again the corresponding remainder \e{eq:Grr1} is in $\ell^1 ({\Bbb Z}_{+})$. Therefore repeating the arguments of Sect.~2 and 3, we find that equation \e{eq:Jy} has a solution 
(the Jost solution) $f_{n} (z)$     distinguished by the asymptotics 
\[
f_{n} (z)=a_{n}^{-1/2} (\sgn\beta_\infty)^n e^{-  \varphi_{n} - z \psi_{n}}  (1+ o(1)), \q   n\to \infty,
\]
%label{eq:D1+}\end{equation} 
where   $\varphi_{n}$ and $\psi_{n}$ are  defined by formulas \e{eq:Grf1} and \e{eq:psi}, respectively.
Now the Jost solution $f_{n} (z)$ is analytic in the whole complex plane so that the spectrum of $J$ is discrete. This is of course consistent with Proposition~\ref{discr}. 
The second solution  $g_{n} (z)$ of the Jacobi equation \e{eq:Jy} is again given by equality \e{eq:GEg} which leads to the asymptotic formula
  \begin{equation}
  P_{n}(z)=-      \Omega(z)  \frac{ (\sgn \beta_{\infty})^{n +1} }{2\sqrt{\beta^2_{\infty}-1}}  \frac{e^{\varphi_{n}+ z\psi_{n}}  }{  \sqrt{a_{n}}} (1+ o(1)), \q n\to\infty.
\label{eq:Gas+}\end{equation}
 As in Sect.~4.2,
the resolvent  of $J$ is determined by formula \e{eq:RRe}, but, in contrast to the case $|\beta_{\infty}|<1$,  its singularities are due to zeros of the denominator $\Omega (z)$ only.
    
           \subsection{The non-Carleman case versus Carleman one}
           
 Let us now compare the results of this paper for the non-Carleman case \e{eq:nc} with those described in the previous subsection.  
 
  Suppose first that $|\beta_{\infty}|< 1$. According to Theorem~\ref{GSS}  all solutions of equation \e{eq:Jy} for every $z\in{\Bbb C}$ belong to $\ell^2 ({\Bbb Z}_{+})$ so that the Jacobi operator $J_{0}$ has a one parameter family of self-adjoint extensions $J$. In this case formula \e{eq:Jy} for the resolvents of the operators $J$ makes no sense. Nevertheless    the scalar products $((J-z)^{-1}e_{0}, e_{0})$, that is,   the Cauchy-Stieltjes transforms of the corresponding spectral measures $d\rho_{J}(\lambda)$,   can be expressed via the orthogonal polynomials of the first  $P_{n}(z)$  and of the second $\wt{P}_{n} (z)$ kinds
   by the Nevanlinna formula   obtained by him in \cite{Nevan} (see, e.g., formulas (2.31), (2.32) in the book \cite{AKH}).
   This remarkable formula implies, in particular, that the spectra of all self-adjoint extensions $  J$ of $J_{0}$ are discrete. The spectral measure at eigenvalues of $J$ is given by formula  \e{eq:meas}. Our construction is quite independent of the Nevanlinna theory, but if some link with this theory exists, it would be desirable to find it.

   Asymptotic formulas \e{eq:A22r} and \e{eq:Sz} look rather similar although
    the phase   in \e{eq:Sz} contains an additional term $\lambda\psi_{n}$.  The amplitude factor $a_{n}^{-1/2}$ in \e{eq:A22r} and \e{eq:Sz} is the same.
However, under assumption \e{eq:Carl} the sequence \e{eq:Sz} never belongs to $\ell^2({\Bbb Z}_{+})$ while under assumption \e{eq:nc} the sequence \e{eq:A22r} belongs to this space for all $\lambda\in{\Bbb R}$.

The difference between the Carleman and non-Carleman cases is even more obvious for $\Im z\neq 0$. According to Theorem~\ref{LC} in the  non-Carleman case,  all solutions of the Jacobi equation oscillate as $n\to\infty$, but they are in $\ell^2 ({\Bbb Z}_{+})$ due to the amplitude factor $a_{n}^{-1/2}$. In the  Carleman case,   the solution $f_{n}  (z)$ exponentially decays while  the solution  $g_{n}  (z)$ exponentially grows as $n\to\infty$.

In the case $|\beta_{\infty}|> 1$, the asymptotic formulas \e{eq:GEGE} and \e{eq:Gas+} for orthogonal polynomials  look 
similar, but the first of them contains an additional factor $\varkappa_{\infty} $ while \e{eq:Gas+} contains a factor $e^{z\psi_{n}}$. We  emphasize that in the Carleman   case \e{eq:Carl} the Jacobi operators $J_{0}$ are essentially self-adjoint for all sequences $b_{n}$ while in the non-Carleman   case \e{eq:nc} the Jacobi operators $J_{0}$ are essentially self-adjoint 
if and only if condition \e{eq:nc+}  is satisfied (see Corollary~\ref{GEc}). This answer is illustrated by Examples~\ref{GEd}  and \ref{GEex}.

  \end{document}